\documentclass[%
 aip,
 jmp,%
 amsmath,amssymb,
 reprint,%
]{revtex4-1}

\usepackage{relsize}
\usepackage{verbatim}
\usepackage{graphicx}
\graphicspath{{Figures_pngNew/}}
\usepackage{dcolumn}
\usepackage{bm}
\usepackage{subfigure}

\usepackage{amsmath,amssymb,amsfonts}
\usepackage{graphics}
\usepackage{color}
\usepackage{xcolor}
\definecolor{rev}{rgb}{0,0,0}

\usepackage[ruled,vlined]{algorithm2e}
\begin{document}

\title{Forward sensitivity approach for estimating eddy viscosity closures in nonlinear model reduction}

\author{Shady E Ahmed}
\author{Kinjal Bhar}
\affiliation{ 
School of Mechanical \& Aerospace Engineering, Oklahoma State University, Stillwater, OK 74078, USA.
}%

\author{Omer San}%
 \email{osan@okstate.edu}
\affiliation{ 
School of Mechanical \& Aerospace Engineering, Oklahoma State University, Stillwater, OK 74078, USA.
}%

\author{Adil Rasheed}%
\affiliation{ 
	Department of Engineering Cybernetics, Norwegian University of Science and Technology, N-7465, Trondheim, Norway.
}%



\date{\today}

\begin{abstract}
In this paper, we propose a variational approach to estimate eddy viscosity using forward sensitivity method (FSM) for closure modeling in nonlinear reduced order models. FSM is a data assimilation technique that blends model's predictions with noisy observations to correct initial state and/or model parameters. We apply this approach on a projection based reduced order model (ROM) of the one-dimensional viscous Burgers equation with a square wave defining a moving shock, \textcolor{rev}{and the two-dimensional vorticity transport equation formulating a decay of Kraichnan turbulence.} We investigate the capability of the approach to approximate an optimal value for eddy viscosity with different measurement configurations. Specifically, we show that our approach can sufficiently assimilate information either through full field or sparse noisy measurements to estimate eddy viscosity closure to cure standard Galerkin reduced order model (GROM) predictions. Therefore, our approach provides a modular framework to correct forecasting error from a sparse observational network on a latent space. We highlight that the proposed GROM-FSM framework is promising for emerging digital twin applications, where real-time sensor measurements can be used to update and optimize surrogate model's parameters.  
\end{abstract}


\keywords{Forward sensitivity method, Galerkin projection, proper orthogonal decomposition, reduced order modeling, closure modeling, moving shock, Burgers equation} 
\maketitle


\section{Introduction} \label{sec:intro}

Data assimilation (DA) is a family of algorithms and techniques that aim at blending mathematical models with (noisy) observations to provide better predictions by correcting initial condition and/or model's parameters \cite{ghil1991data,kalnay2003atmospheric,lewis2006dynamic}. DA plays a key role in geophysical and meteorological sciences to make more reliable numerical weather predictions. Standard popular algorithms that are often adopted in weather prediction centers include variational methods (e.g., 3D-VAR \cite{lorenc1986analysis,parrish1992national} and 4D-VAR \cite{courtier1997dual,rabier2000ecmwf,elbern20004d,courtier1994strategy,lorenc2005does,gauthier2007extension} methods), sequential methods (e.g., reduced rank (ensemble) Kalman filters \cite{houtekamer1998data,burgers1998analysis,evensen2003ensemble,houtekamer2001sequential,houtekamer2005ensemble,treebushny2003new,buehner2003reduced,lakshmivarahan2009ensemble}), and hybrid methods \cite{zupanski2005maximum,desroziers20144denvar,lorenc2015comparison,wang2008hybrid,buehner2013four,kleist2015osse,kleist2015osse2}. Another method that mitigates the computational cost in solving the inherent optimization problem in variational methods is called the forward sensitivity method (FSM) developed by Lakshmivarahan and Lewis \cite{lakshmivarahan2010forward,lakshmivarahan2017forecast}. FSM builds on the assumption that model error stems from incorrect specification of the control elements, which include initial conditions, boundary conditions, and physical/empirical parameters. The FSM approach corrects the control elements using information from the time evolution of sensitivity functions, defined as the derivatives of model output with respect to the elements of control.

Other than meteorology \cite{wang2000data}, DA tools are gaining popularity in different disciplines like reservoir engineering \cite{aanonsen2009ensemble}, and neuroscience \cite{hutt2018forecast}. Recent works have also drawn techniques and ideas from DA to enrich reduced order modeling of fluid flows and vice versa \cite{protas2015optimal,zerfas2019continuous,xiao2018parameterised,daescu2007efficiency,cstefuanescu2015pod,cao2007reduced,robert2005reduced,arcucci2019optimal,popov2020multifidelity}. In conventional projection-based model reduction approaches, a set of system's realizations are used to build a reduced order model (ROM) that sufficiently represent the system's dynamics with significantly lower computational cost \cite{bai2002krylov,lucia2004reduced,hess2019localized,kramer2019nonlinear,swischuk2019projection,bouvrie2017kernel,hamzi2019local,korda2018data,korda2018linear,hartmann2018model,holmes2012turbulence,taira2017modal,taira2019modal,noack2011reduced,rowley2017model,nair2019transported,kaiser2014cluster,haasdonk2011training}. This process includes the extraction of a handful of basis functions representing the underlying flow patterns or coherent structures that dominate the majority of the bulk mass, momentum and energy transfers. \textcolor{rev}{ In the fluid dynamics research community}, proper orthogonal decomposition (POD) is, generally speaking, the most popular and effective technique that produces hierarchically ordered solution-adapted basis functions (or modes) that provide the optimal basis to represent a given collection of field data or snapshots \cite{sirovich1987turbulence,berkooz1993proper,holmes2012turbulence,chatterjee2000introduction,rathinam2003new}. To emulate \textcolor{rev}{the system's dynamics}, a surrogate model is often built by performing a Galerkin projection of the full order model (FOM) operators onto a reduced subspace spanned by the formerly constructed POD modes \cite{ito1998reduced,iollo2000stability,rowley2004model,milk2016pymor,puzyrev2019pyrom,bergmann2009enablers,couplet2005calibrated,kunisch2001galerkin}. 

However, the off-design performance of ROMs is usually questionable since the reduced basis and operators are formed offline for a given set of operating conditions, while the ROM has to be solved online for different conditions. Therefore, a dynamic update of model operators and parameters is often sought to enhance the applicability of ROMs in realistic contexts. That being said, adoption of DA tools to absorb real observations to correct and update ROMs should present a viable cure for this caveat. 
The present paper aims at pushing towards utilizing DA techniques to improve the performance of nonlinear ROMs. \textcolor{rev}{A common problem that emerges in such ROMs is the inaccuracy of solution trajectory, especially for long time integration of quasi-stationary problems. This solution inaccuracy has been commonly attributed to the modal truncation and intrinsic interactions between truncated modes and retained modes. A correction term compensating the effects of truncated modes has been often introduced to achieve more accurate ROM results \cite{wang2012proper,san2014proper,osth2014need,baiges2015reduced,kondrashov2015data,fick2018stabilized}. Furthermore, recent studies have shown that model's performance can be improved by the choice of the projection method \cite{oberai2016approximate,choi2019space,grimberg2020stability} and the definition of the adopted inner product \cite{kalashnikova2010stability}.}

In order to enhance the solution accuracy, closure and stabilization techniques have been introduced to account for the effects of discarded modes on the dynamics of the ROM. In particular, eddy viscosity closures, inspired from large eddy simulations (LES), have shown a significant success in ROM closure modeling \cite{borggaard2011artificial,wang2011two,akhtar2012new,cordier2013identification,protas2015optimal,san2015stabilized}. The estimation of an optimal value of the eddy viscosity parameter has been the topic for many research works though. For example, empirical relations can be adopted \cite{rempfer1997proper,san2014proper,ahmed2018stabilized}, or ideas can be borrowed from LES frameworks to dynamically compute a better approximation of the eddy viscosity parameter \cite{wang2012proper,rahman2019dynamic,mou2020data,imtiaz2020nonlinear}. Moreover, a 4D-VAR approach has been suggested to provide an optimal nonlinear eddy viscosity estimate in Galerkin projection based ROMs\cite{protas2015optimal}. An adaptive nudging technique has also been recently introduced to force ROMs towards the reference solution corresponding to the observed data \cite{zerfas2019continuous}. 

Instead, in the present paper, we propose a novel framework to estimate eddy viscosity closure using noisy observations from a sparse observation network. In particular, we adopt the forward sensitivity method to evaluate the sensitivity of ROM predictions to the eddy viscosity parameter. Observations, whenever available, can therefore be used to approximate a more representative value of eddy viscosity to better reflect the true system's dynamics. We highlight that the proposed approach is very suitable for emerging digital twin applications \cite{tao2018digital,tao2018digitalb,madni2019leveraging,rasheed2020digital,ganguli2020digital,chakraborty2020role}, where real-time measurements are abundant (and noisy). Thus, efficiently assimilating these measurements to improve ROMs can be a key enabler for such applications which require many-query and near real-time simulations. We test our approach using \textcolor{rev}{two test cases of varying complexity, namely the one-dimensional viscous Burgers equation with a square wave representing a moving shock and the two-dimensional vorticity transport equation applied to Kraichnan turbulence.} We apply the proposed GROM-FSM to assimilate information from either full field or sparse field measurements. Therefore, our approach provides a modular framework to optimally estimate closure parameters for submodal scale physics, which can be effectively used in emerging sensor-centric applications in transport processes.

The rest of the paper is outlined here. In Section~\ref{sec:fsm}, we review the forward sensitivity method and its mathematical foundation as an established data assimilation algorithm. We then construct the standard Galerkin ROM and the corresponding reduced operators for the 1D Burgers problem \textcolor{rev}{and the 2D vorticity transport equation} in Section~\ref{sec:rom}. Then, we describe the proposed approach for closure estimation via FSM, namely GROM-FSM, in Section~\ref{sec:closure}. Results and relevant discussions are provided in Section~\ref{sec:results}. In particular, we consider the assimilation of full field and sparse field measurements. For the latter, we explore two approaches for assimilating information from sparse observations. \textcolor{rev}{We also extend the eddy viscosity estimation framework to permit mode-dependent closures.} Concluding remarks and insights are drawn in Section~\ref{sec:conc}.

\section{Forward Sensitivity Method} \label{sec:fsm}
In this section, we briefly describe the forward sensitivity method (FSM) proposed by Lakshmivarahan and Lewis \cite{lakshmivarahan2010forward}. The idea behind this technique is to find optimal control parameters by iteratively correcting the control for the least squares fit of the model to the observational data. The control parameters in question here can be any unknown such as initial conditions, boundary conditions, and physical model parameters. The correction to each control parameter is dictated by its corresponding sensitivity function. In essence, the sensitivity function is the quantitative measure of influence of each control parameter on the model states. \textcolor{rev}{The nature} of combining physical models with actual data to solve an inverse problem is what makes FSM a modular DA approach. 

Let the dynamical system of interest be defined by a set of ordinary differential equations (ODEs) as below,
\begin{equation} \label{eq:dyn_sys}
     \dfrac{d\mathbf{x}}{dt} = \mathbf{f}(\mathbf{x},\boldsymbol \alpha),
\end{equation}
where $\mathbf{x}(t) \in \mathbb{R}^{n}$  is the system state-vector with the initial condition $\mathbf{x}^0$ and $\mathbf{\boldsymbol \alpha} \in \mathbb{R}^{p}$ denotes the physical parameters. The vector of control parameters is represented as $\mathbf{c}=[\mathbf{x}^0,\mathbf{\boldsymbol \alpha}]^{T} \in \mathbb{R}^{n+p}$. Here, it is assumed that the solution $\mathbf{x}(t)$ exists and is unique and has a smooth dependence with the control vector $\mathbf{c}$.

Discretizing Eq.~\ref{eq:dyn_sys} by using some numerical method like Runge-Kutta schemes, we get a model equation which gives the evolution of model states in discrete time as,
\begin{equation} \label{eq:dis_dyn_eqn}
    \mathbf{x}^{k+1} = \mathbf{M}(\mathbf{x}^k,\boldsymbol \alpha),
\end{equation}
where $\mathbf{x}^k = [x_{1}^{k}, x_{2}^{k}, \dots, x_{n}^{k}]^{T}$ denotes the time-discretized model states at discrete time $t_k$ and $\mathbf{M}=[M_1(\mathbf{x}^k,\boldsymbol \alpha), M_2(\mathbf{x}^k,\boldsymbol \alpha), \dots, M_n(\mathbf{x}^k,\boldsymbol \alpha)]^{T}$ refer to the state transition maps from time $t_k$ to $t_{k+1}$. Differentiating Eq.~\ref{eq:dis_dyn_eqn} with respect to $\mathbf{x}^0$, we get 
\begin{equation} \label{eq:IC_sen}
    \dfrac{\partial x^{k+1}_i}{\partial x^{0}_j} = \sum^{n}_{q=1} \left( \dfrac{\partial M_i}{\partial x^k_q} \right) \left( \frac{\partial x^k_q}{\partial x^0_j} \right),
\end{equation}
where $1 \leqslant i,j \leqslant n$. Similarly, differentiating Eq.~\ref{eq:dis_dyn_eqn} with respect to $\mathbf{\boldsymbol \alpha}$, we obtain
\begin{equation} \label{eq:para_sen}
    \frac{\partial x^{k+1}_i}{\partial \alpha_j} = \sum^{n}_{q=1} \left( \frac{\partial M_i}{\partial x^k_q} \right) \left( \frac{\partial x^k_q}{\partial \alpha_j} \right) + \frac{\partial M_i}{\partial \alpha_j} 
\end{equation}
where $1 \leqslant i \leqslant n$ and $1 \leqslant j \leqslant p$. In Eq.~\ref{eq:IC_sen} and Eq.~\ref{eq:para_sen}, the superscript refers to the discrete time index while the subscript refers to the specific component. Now, we can define $\mathbf{U}^k$ as the sensitivity matrix of $\mathbf{x}^k$ with respect to initial state, where $[\mathbf{U}^k]_{ij} = \partial x^k_i/\partial x^0_j$ for $1 \leqslant i,j \leqslant n$. Also, we define $\mathbf{V}^k$ as the sensitivity matrix of $\mathbf{x}^k$ with respect to the parameter-vector $\mathbf{\boldsymbol \alpha}$, where $ [\mathbf{V}^k]_{ij}  = \partial x^k_i/\partial \alpha_j$ for $1 \leqslant i \leqslant n$ and $1 \leqslant j \leqslant p$. Then, we can rewrite Eqs.~\ref{eq:IC_sen}-~\ref{eq:para_sen} in matrix as below
\begin{align} 
    \mathbf{U}^{k+1} &= \mathbf{D}_\mathbf{x}^k \mathbf{(M)}\mathbf{U}^{k}, \label{eq:evolutio_Uk} \\
    \mathbf{V}^{k+1} &= \mathbf{D}_\mathbf{x}^k \mathbf{(M)}\mathbf{V}^{k} + \mathbf{D_{\boldsymbol \alpha}(M)}, \label{eq:evolutio_Vk}
\end{align}
initialized as $\mathbf{U}^0 = \mathbf{I}$ and $\mathbf{V}^0 = \mathbf{0}$.

Here, $\mathbf{D}_\mathbf{x}^k \mathbf{(M)}$ and $\mathbf{D_{\boldsymbol \alpha}(M)}^k$ are the Jacobian matrices of $\mathbf{M(\cdot)}$ with respect to $\mathbf{x}$ and $\boldsymbol \alpha$ at discrete time $t_k$, respectively. 
Moreover, $\mathbf{U}^k \in \mathbb{R}^{n \times n}$ and $\mathbf{V}^k  \in \mathbb{R}^{n \times p}$ are called the forward sensitivity matrices with respect to initial conditions and parameters, respectively. In effect, the system dynamics in Eq.~\ref{eq:dis_dyn_eqn} gets reduced to a set of linear matrix equations (Eq.~\ref{eq:evolutio_Uk} and Eq.~\ref{eq:evolutio_Vk}) which give the evolution of the sensitivity matrices in discrete time. By first order approximation, we have
\begin{equation} \label{eq:dx}
    \Delta \mathbf{x}^k \approx \delta \mathbf{x}^k = \mathbf{U}^k \delta \mathbf{x}^0 + \mathbf{V}^k \delta \boldsymbol \alpha , 
\end{equation}
where $\delta\mathbf{x} \in \mathbb{R}^{n}$. 

So far, no observational data have been used. Let $\mathbf{z}(t) \in \mathbb{R}^{m}$ be the observation vector available for $N$ time snapshots; and $\mathbf{h}: \mathbb{R}^n \rightarrow \mathbb{R}^m$ maps the model space $\mathbb{R}^{n}$ to the observation space $\mathbb{R}^{m}$. Hence, the observation vector can be defined mathematically as follows,
\begin{equation} \label{eq:obs map}
     \textcolor{rev}{\mathbf{z}(t) = \mathbf{h}(\Tilde{\mathbf{x}}(t)) + \mathbf{v}(t)},
\end{equation}
where $\textcolor{rev}{\Tilde{\mathbf{x}}(t)} \in \mathbb{R}^{n}$ is the true state of the system and $\mathbf{v}(t) \in \mathbb{R}^{m}$ represents the measurement noise, which is assumed to be white Guassian noise with zero mean and covariance matrix $\mathbf{R}(t) \in \mathbb{R}^{m \times m}$. Writing Eq.~\ref{eq:obs map} in the discrete-time form we get, 
\begin{equation} \label{eq:dis_obs_map}
    \mathbf{z}^k = \mathbf{h}(\Tilde{\mathbf{x}}^k) + \mathbf{v}^k,
\end{equation}
where $\mathbf{v}^k$ is white Gaussian noise with the covariance matrix $\mathbf{R}^k$. In most cases, $\mathbf{R}^k$ is a diagonal matrix. For simplicity, we assume that $\mathbf{R}^k= \sigma_{Obs}^2 \mathbf{I}_m$, where $\mathbf{I}_m$ is the $m\times m$ identity matrix. 

Assuming that the model is \textcolor{rev}{a perfect representation} of the actual physical phenomenon and given a starting guess value of the control $\mathbf{c}$, we can run the model forward to predict $\mathbf{x}^k$ $\forall$ $1 \leqslant k \leqslant N$, then the forecast error $\mathbf{e}_F^k \in \mathbb{R}^{m}$ defined as,
\begin{equation} \label{eq:error}
    \mathbf{e}_F^k = \mathbf{z}^k - \mathbf{h}(\mathbf{x}^k).
\end{equation}

The forecast error $\mathbf{e}_F^k$ is composed of the sum of a deterministic part defined as $\mathbf{h}(\Tilde{\mathbf{x}}^k) - \mathbf{h}(\mathbf{x}^k)$ and a random part $\mathbf{v}^k$. The random error stems from the inherent error in the mapping $\mathbf{h}: \mathbf{x}^k \rightarrow \mathbf{z}^k$ and we have no control on it, however it is the goal of FSM to minimize the deterministic part in a least squares sense at all the $N$ time snaps by choosing an optimal value for $\mathbf{c}$. 

Now, the goal of FSM is to find a perturbation to the control $\delta \mathbf{c}$ from the given starting guess $\mathbf{c}$. This, in turn, would cause a $\delta \mathbf{x}^k$ change in $\mathbf{x}^k$ such that the actual observation matches with the forecast observation from the model as follows,
\begin{equation}
    \mathbf{z}^k = \mathbf{h}(\mathbf{x}^k + \delta \mathbf{x}^k) \approx \mathbf{h}(\mathbf{x}^k) + \mathbf{D}^k_{\mathbf{x}}(\mathbf{h}) \delta \mathbf{x}^k.
\end{equation}
Thus, the forecast error $\mathbf{e}_F^k$ can be written as, 
\begin{equation} \label{eq:e_dx}
    \mathbf{e}_F^k = \mathbf{D}^k_{\mathbf{x}}(\mathbf{h}) \delta \mathbf{x}^k.
\end{equation}

Combining Eq.~\ref{eq:dx} with Eq.~\ref{eq:e_dx}, and setting $\mathbf{H}_1^k = \mathbf{D}^k_{\mathbf{x}}(\mathbf{h}) \mathbf{U}^k \in \mathbb{R}^{m \times n}$, $\mathbf{H}_2^k = \mathbf{D}^k_{\mathbf{x}}(\mathbf{h}) \mathbf{V}^k \in \mathbb{R}^{m \times p}$, we get, 
\begin{equation} \label{eq:e_dx_da}
     \mathbf{e}_F^k = \mathbf{H}_1^k \delta \mathbf{x}^0  + \mathbf{H}_2^k \delta \boldsymbol{\alpha}.
\end{equation}
Equation~\ref{eq:e_dx_da} can be further simplified and written in terms of the perturbation to the control $\delta \mathbf{c}$ as
\begin{equation} \label{eq:e_dc}
    \mathbf{H}^k \delta \mathbf{c} = \mathbf{e}_F^k,
\end{equation}
where $\mathbf{H}^k = [ \mathbf{H}_1^k,\mathbf{H}_2^k] \in \mathbb{R}^{m \times (n+p)}$ and $\delta \mathbf{c} = [ \delta \mathbf{x}^0 , \delta \boldsymbol{\alpha} ]^T \in \mathbb{R}^{n+p}$.

Equation~\ref{eq:e_dc} can be formulated for all the $N$ time snapshots for which observations are available and the following linear equation is obtained,
\begin{equation} \label{eq:E_dc}
    \mathbf{H} \delta \mathbf{c} = \mathbf{e}_F,
\end{equation}
where the matrix $\mathbf{H} \in \mathbb{R}^{N m \times (n+p)}$ and the vector $\mathbf{e}_F \in \mathbb{R}^{N m}$ are defined as follows,
\begin{equation} \label{eq:H_e_def}
    \mathbf{H} = \begin{bmatrix}     
    \mathbf{H}^1 \\
    \mathbf{H}^2 \\
    \vdots \\
    \mathbf{H}^{N}   
    \end{bmatrix}, \qquad
    \mathbf{e}_F = \begin{bmatrix}
    \mathbf{e}_F^1 \\
    \mathbf{e}_F^2 \\
    \vdots \\
    \mathbf{e}_F^{N}
    \end{bmatrix}.
\end{equation}
Depending on the value of $Nm$ relative to $(n+p)$, Eq.~\ref{eq:E_dc} can give rise to either an over-determined or an under-determined linear inverse problem. In either case, the inverse problem can be solved in a weighted least squares sense to find an optimal value of $\delta \mathbf{c}$, with $\mathbf{R}^{-1}$ as a weighting matrix, where $\mathbf{R}$ is a block-diagonal matrix constructed as follows,
\begin{equation}
    \mathbf{R} = \begin{bmatrix}
                 \mathbf{R}^{1} &              &        & \\
                                & \mathbf{R}^2 &        & \\
                                &              & \ddots & \\
                                &              &        & \mathbf{R}^N
                \end{bmatrix}.
\end{equation}

For simplicity, we assume that $\mathbf{R}$ is a diagonal matrix defined as $\mathbf{R} = \sigma_{Obs}^2 \mathbf{I}_{Nm}$, where $\mathbf{I}_{Nm}$ is the $Nm \times Nm$ identity matrix. Then, the solution of Eq.~\ref{eq:E_dc} can be written as
\begin{equation}\label{eq:LSsolve}
\delta \mathbf{c} = \begin{cases}
                     \left( \mathbf{H}^T \mathbf{R}^{-1} \mathbf{H} \right)^{-1} \mathbf{H}^T \mathbf{R}^{-1} \mathbf{e}_{F}, &\quad\text{over-determined,} \\
                    \mathbf{R}^{-1} \mathbf{H}^T  \left( \mathbf{H} \mathbf{R}^{-1} \mathbf{H}^T \right)^{-1} \mathbf{e}_{F}, &\quad\text{under-determined.}
                    \end{cases}
\end{equation}
It has been seen that the first order approximation progressively yield better results by repeating the entire process for multiple iterations until convergence with certain tolerance \cite{lakshmivarahan2010forward}. 



\section{Reduced Order Modeling} \label{sec:rom}
\textcolor{rev}{In this section, we briefly derive a reduced order model (ROM) for a dynamical system governed by the following autonomous partial differential equation (PDE)
\begin{equation} \label{eq:pde}
    \dfrac{\partial q}{\partial t} = \mathcal{F}(q),
\end{equation}
where $q$ is the state of system (flow field variables) and $\mathcal{F}(q)$ governs the dynamics of $q$.}
We follow the standard Galerkin projection to construct the sought ROM which includes two main steps. First, the flow field variable $q(\mathbf{x},t)$ (where $q(\mathbf{x},t)$ represents the vectorized form of $q$ at time $t$) is approximated as a linear superposition of the contributions of a few modes, which can be mathematically expressed as
\begin{equation} \label{eq:qROM}
    q(\mathbf{x},t) = \bar{q}(\mathbf{x}) + \sum_{k=1}^{R} a_k(t) \phi_k(\mathbf{x}),
\end{equation}
where $\bar{q}(\mathbf{x})$ represents the mean-field, $\phi_k(\mathbf{x})$ are the spatial modes (or basis functions), $a_k(t)$ are the time-dependent modal coefficients (i.e., weighting functions), and $R$ is the number of retained modes in ROM approximation (i.e., ROM dimension). The second step is to project the governing equation (i.e., Eq.~\ref{eq:pde}) onto the subspace spanned by $\{\phi_k\}_{k=1}^R$. Thus, the two main ingredients for building a \textcolor{rev}{Galerkin ROM} (GROM) are the basis functions $\{\phi_k\}_{k=1}^R$ and a Galerkin projection of the governing equation. To compute the basis functions $\{\phi_k\}_{k=1}^R$, we follow the popular proper orthogonal decomposition (POD) approach described in Section~\ref{sec:POD}, followed by derivation of GROM equations in Section~\ref{sec:GP}. \textcolor{rev}{Overall, this GROM approach utilizes a linear decomposition technique that is able to properly treat the nonlinearity of $\mathcal{F}(q)$, since it accounts for nonlinear coupling of terms acting within the linear space defined by the POD basis functions \cite{lucia2004reduced}.}

\subsection{Proper Orthogonal Decomposition} \label{sec:POD}
Proper orthogonal decomposition (POD) is \textcolor{rev}{a data-driven modal decomposition technique that gained remarkable popularity in the fluid mechanics community} due to its simplicity as well as robustness. Given a set of solution trajectories or realizations (known as snapshots), POD lays out a systematic approach to compute a solution-adapted basis functions that provide the optimal basis to represent a given set of simulation data or snapshots. \textcolor{rev}{Specifically, POD produces hierarchically organized basis functions}, based on their contribution to the total system's energy, which makes the modal selection a trivial process. In particular, given a collection of system realizations, we build a snapshot matrix $\mathbf{A} \in \mathbb{R}^{n \times N}$ as follows,
\begin{equation}
 \mathbf{A} = \begin{bmatrix}
 q(x_1,t_1) & q(x_1,t_2) & \dots & q(x_1,t_{N}) \\
 q(x_2,t_1) & q(x_2,t_2) & \dots & q(x_2,t_{N}) \\
 \vdots     & \vdots    &  \ddots      & \vdots \\
 q(x_{n},t_1) & q(x_{n},t_2) & \dots & q(x_{n},t_{N}) \\
 \end{bmatrix},
\end{equation}
where $n$ is the number of spatial locations and $N$ is the number of snapshots. \textcolor{rev}{A mean-subtracted snapshot matrix $\tilde{\mathbf{A}}$ is defined as $\tilde{\mathbf{A}} = \mathbf{A} - \dfrac{1}{N} \mathbf{A} \mathbf{1}_{N\times N}$, where $\mathbf{1}_{N\times N}$ is an ${N\times N}$ matrix of ones. Then, a thin singular value decomposition (SVD) is performed on $\tilde{\mathbf{A}}$ as follows} , 
\begin{equation}
    \tilde{\mathbf{A}} = \mathbf{U} \mathbf{\Sigma} \mathbf{V}^T,
\end{equation}
where $\mathbf{U} \in \mathbb{R}^{n \times N}$ is a matrix with orthonormal columns are the left singular vectors of \textcolor{rev}{$\tilde{\mathbf{A}}$}, which represent the spatial basis as,
\begin{equation}
 \mathbf{U} = \begin{bmatrix}
 U_1(x_1) & U_2(x_1) & \dots & U_{N}(x_1) \\
 U_1(x_2) & U_2(x_2) & \dots & U_{N}(x_2) \\
 \vdots     & \vdots    &  \ddots      & \vdots \\
 U_1(x_{n}) & U_2(x_{n}) & \dots & U_{N}(x_{n})
 \end{bmatrix},
\end{equation}
while the columns of $\mathbf{V} \in \mathbb{R}^{N \times N}$ are the right singular vectors of \textcolor{rev}{$\tilde{\mathbf{A}}$}, representing the temporal basis as
\begin{equation}
 \mathbf{V} = \begin{bmatrix}
 V_1(t_1) & V_2(t_1) & \dots & V_{N}(t_1) \\
 V_1(t_2) & V_2(t_2) & \dots & V_{N}(t_2) \\
 \vdots     & \vdots    &  \ddots      & \vdots \\
 V_1(t_{N}) & V_2(t_{N}) & \dots & V_{N}(t_{N})
 \end{bmatrix}.
\end{equation}
The singular values of \textcolor{rev}{$\tilde{\mathbf{A}}$} are stored in descending order as the entries of the diagonal matrix $\mathbf{\Sigma} \in \mathbb{R}^{N \times N}$,
\begin{equation}\label{eq:sigma}
 \mathbf{\Sigma} = \begin{bmatrix}
 \sigma_1 &          &        &  \\
          & \sigma_2 &        &  \\
          &          & \ddots &  \\
          &          &        & \sigma_{N}
 \end{bmatrix},
\end{equation}
where $\sigma_1 \ge \sigma_2 \ge \dots \sigma_{N} \ge 0$. For dimensionality reduction purposes, only the first $R$ columns of $\mathbf{U}$, corresponding to the largest $R$ singular values, are stored. Those represent the most effective $R$ POD modes, denoted as $\{\phi_k\}_{k=1}^{R}$ in the rest of the manuscript. The computed basis functions are orthonormal by construction as
\begin{equation}
\langle \phi_i ; \phi_j \rangle = 
     \begin{cases}
       1 &\quad\text{if } i = j \\
       0 &\quad\text{otherwise,}
     \end{cases}    
\end{equation}
where the angle parentheses $\langle \cdot ; \cdot \rangle$ stands for the standard inner product in Euclidean space (i.e., dot product). We note that the presented direct algorithm might be unfeasible for larger data sets, as stacking snapshots into a single huge matrix is usually prohibitive. Instead, the method of snapshots \cite{sirovich1987turbulence} can be followed to efficiently approximate the POD bases.

\subsection{Galerkin ROM} \label{sec:GP}
Having a set of POD basis functions in hand, an orthogonal projection can be performed to obtain the \textcolor{rev}{Galerkin ROM} (GROM). To do so, the ROM approximation (Eq.~\ref{eq:qROM}) is substituted into the governing equation and an inner product with the POD basis functions is carried out. In deriving the GROM equations, we highlight that the POD bases are only spatial functions (i.e., independent of time) and the modal coefficients are independent of space. We also utilize the orthonormality property of the basis functions to get the following set of ordinary differential equations (ODEs) representing the tensorial GROM
\textcolor{rev}{\begin{align}
    \dfrac{\mathrm{d}a_k}{\mathrm{d}t} &=   \mathfrak{B}_k + \sum_{i=1}^{R} \mathfrak{L}_{i,k} a_i + \sum_{i=1}^{R} \sum_{j=1}^{R} \mathfrak{N}_{i,j,k} a_i a_j, \label{eq:grom}
\end{align}
where $\mathfrak{B}$, $\mathfrak{L}$ and $\mathfrak{N}$ are the vector, matrix and tensor of predetermined model coefficients corresponding to constant, linear and nonlinear terms, respectively. We note here that the last term results from the quadratic nonlinearity encountered in most of the fluid flow systems. In particular, we consider here two cases of particular interest. First, we consider the one-dimensional Burgers equation as a prototypical test bed for transport systems with quadratic nonlinearity and Laplacian dissipation. For this case, we solve the problem of a moving shock, which can be considered as a challenging case for ROM applications \cite{ahmed2018stabilized}. In the second case, we consider the vorticity transport equation, with an application to the two-dimensional decaying turbulence. } 

\textcolor{rev}{
\subsubsection{1D Burgers problem} \label{sec:1Dbrg}
The one-dimensional (1D) Burgers equation is defined with the following partial differential equation (PDE)
\begin{equation} \label{eq:burgers}
     \dfrac{\partial u}{\partial t}  = - u \dfrac{\partial u}{\partial x} + \dfrac{1}{\text{Re}} \dfrac{\partial^2 u}{\partial x^2},
\end{equation}
where $\text{Re}$ is Reynolds number relating the inertial and viscous effects. Using the following definition
\begin{equation} \label{eq:burgers2}
      u(x,t) = \bar{u}(x) + \sum_{k=1}^R a_k(t) \phi_k(x),
\end{equation}
the GROM model coefficients can be precomputed during an offline stage as
\begin{align*}
  \mathfrak{B}_{k} & = \big\langle -\bar{u} \dfrac{\partial \bar{u}}{\partial x} + \dfrac{1}{\text{Re}} \dfrac{\partial^2 \bar{u} }{\partial x^2} ;  \phi_k \big\rangle, \\
    \mathfrak{L}_{i,k} & = \big\langle -\bar{u} \dfrac{\partial \phi_i}{\partial x} - \phi_i \dfrac{\partial \bar{u}}{\partial x} + \dfrac{1}{\text{Re}}\dfrac{\partial^2 \phi_i }{\partial x^2}  ;  \phi_k \big\rangle, \\
    \mathfrak{N}_{i,j,k} &= \big\langle -\phi_i \dfrac{\partial \phi_j}{\partial x};  \phi_k \big\rangle.
\end{align*} }

\textcolor{rev}{
\subsubsection{2D Kraichnan turbulence} \label{sec:2Dturb}
The two-dimensional (2D) Kraichnan turbulence problem models how randomly generated vortices evolve \cite{tabeling2002two}. 
Despite the apparent simplicity, the decaying 2D Kraichnan turbulence is very rich in its dynamics and follows the 2D Navier-Stokes equations, which can be written in vorticity-streamfunction formulation (vorticity-transport equation) as follows
\begin{align} 
\dfrac{\partial \omega}{\partial t}  &= - J(\omega,\psi) + \dfrac{1}{\text{Re}} \nabla^2 \omega, \label{eq:NS}
\end{align}
where $\omega$ is the vorticity and $\psi$ is the streamfunction.  $J(\omega,\psi)$ and $\nabla^2 \omega$ are the Jacobian and Laplacian operators, respectively, which can be defined as
\begin{align}
    J(\omega,\psi) &= \dfrac{\partial \omega}{\partial x} \dfrac{\partial \psi}{\partial y} -  \dfrac{\partial \omega}{\partial y} \dfrac{\partial \psi}{\partial x}, \\
    \nabla^2 \omega &= \dfrac{\partial^2 \omega}{\partial x^2} + \dfrac{\partial^2 \omega}{\partial y^2}.
\end{align}
The vorticity-streamfunction formulation enforces the incompressibility condition, where the vorticity and streamfunction fields are related by the following Poisson equation
\begin{equation}\label{eq:Poisson}
\nabla^2 \psi = -\omega.
\end{equation}}

\textcolor{rev}{With the POD algorithm implemented, a set of POD basis functions $\{\phi_{k}(x,y)\}_{k=1}^{R}$ are obtained from the snapshots of vorticity fields. 
The prognostic variable vorticity in Eq.~\ref{eq:NS} is defined as
\begin{equation} \label{eq:sPOD2}
    \omega(x,y,t) = \bar{\omega}(x,y) + \sum_{k=1}^R a_k(t) \phi_k(x,y). \\
\end{equation}
Since the vorticity and streamfunction are related by the kinematic relationship given by Eq.~\ref{eq:Poisson}, the basis functions ($\theta_k(x,y)$) and mean field ($\bar{\psi}(x,y)$) corresponding to the streamfunction can be obtained from those of the vorticity as follows,
\begin{align}
        &\nabla^2 \bar{\psi}(x,y) = -\bar{\omega}(x,y), \\
        &\nabla^2 \theta_k(x,y) = -\phi_k(x,y), \quad k=1,2,\dots, R,
    \end{align}
which might result in a set of basis functions for the streamfunction that are not necessarily orthonormal. Moreover, the reduced order approximation of streamfunction shares the same temporal coefficients $a_k(t)$,
\begin{equation} \label{eq:sPOD}
    \psi(x,y,t) = \bar{\psi}(x,y) + \sum_{k=1}^R a_k(t) \theta_k(x,y). \\
\end{equation}
}

\textcolor{rev}{The GROM coefficients for this case can be defined as follows \cite{ahmed2020breaking},
\begin{eqnarray}
  & & \mathfrak{B}_{k} = \big\langle  - J(\bar{\omega},\bar{\psi})  + \frac{1}{\mbox{Re}}\nabla^2 \bar{\omega}  ; \phi_{k} \big\rangle , \nonumber \\
  & & \mathfrak{L}_{i,k} = \big\langle - J(\bar{\omega},\theta_{i}) -  J(\phi_{i},\bar{\psi}) + \frac{1}{\mbox{Re}}\nabla^2 \phi_{i}; \phi_{k} \big\rangle  , \nonumber\\
  & &  \mathfrak{N}_{i,j,k} =  \big\langle -J(\phi_{i},\theta_{j}); \phi_{k} \big\rangle. \label{eq:roma7}
\end{eqnarray}}

Due to the quadratic nonlinearity in the aforementioned systems, the computational cost of solving Eq.~\ref{eq:grom} is $O(R^3)$. Therefore, the number of retained modes has to be reduced as much as possible to keep the computational cost affordable. However, this truncation ignores the dyadic interactions between the first $R$ modes and the remaining ones. As a result, an erroneous behaviour might arise in the ROM solution \cite{lassila2014model,rempfer2000low,noack2003hierarchy,gunzburger2019evolve}, and closure/stabilization techniques have been introduced to improve ROM accuracy \cite{balajewicz2012stabilization,wang2012proper,amsallem2012stabilization,kalb2007intrinsic,cordier2013identification}. \textcolor{rev}{As highlighted earlier in Section~\ref{sec:intro}, closure approaches based on physical significance have been usually relied on the analogy between LES and ROMs, e.g., the addition of an artificial viscosity term \cite{borggaard2011artificial}. Recently, data-driven closure methods have been also pursued, e.g., using variational multiscale techniques \cite{stabile2019reduced,reyes2020projection,mou2020data}, machine learning algorithms \cite{pawar2020evolve,san2018extreme,san2018neural,mcquarrie2020data}, and polynomial approximations \cite{xie2018data,mohebujjaman2019physically}.  }

\section{Closure estimation via FSM} \label{sec:closure}
In order to stabilize the GROM, a closure model is usually necessary for complex flows. In the present paper, we consider adding a linear eddy-viscosity term to the governing equation as follows,
\textcolor{rev}{\begin{align} 
    \text{1D Burgers:}& \quad \dfrac{\partial u}{\partial t} = - u \dfrac{\partial u}{\partial x} + (\nu+\nu_e) \dfrac{\partial^2 u}{\partial x^2}, \label{eq:burgersC}\\
    \text{2D Turbulence:}& \quad \dfrac{\partial \omega}{\partial t}  = - J(\omega,\psi) + (\nu+\nu_e) \nabla^2 \omega, \label{eq:NSC}
\end{align}
where $\nu$ is the physical (kinematic) viscosity and $\nu_e$ is an (artificial) eddy viscosity to add an extra dissipation to stabilize the system.} If we follow the same procedure in Section~\ref{sec:GP}, we get the following GROM with closure,
\textcolor{rev}{\begin{align}
    \dfrac{\mathrm{d}a_k}{\mathrm{d}t} &=   \mathfrak{B}_k + \nu_e \widehat{\mathfrak{B}}_k +  \sum_{i=1}^{R} \mathfrak{L}_{i,k} a_i + \nu_e \sum_{i=1}^{R} \widehat{\mathfrak{L}}_{i,k} a_i \nonumber \\
    & + \sum_{i=1}^{R} \sum_{j=1}^{R} \mathfrak{N}_{i,j,k} a_i a_j, \label{eq:gromC}
\end{align}
where $\widehat{\mathfrak{B}}$ and $\widehat{\mathfrak{L}}$ are the constant and linear coefficients resulting from the introduction of the eddy viscosity term and defined as follows
\begin{align*} 
    \text{1D Burgers:}& \quad  \widehat{\mathfrak{B}}_k =  \big\langle \dfrac{\partial^2 \bar{u} }{\partial x^2} ;  \phi_k \big\rangle, \quad \widehat{\mathfrak{L}}_{i,k} = \big\langle  \dfrac{\partial^2 \phi_i }{\partial x^2}  ;  \phi_k \big\rangle,\\
    \text{2D Turbulence:}& \quad  \widehat{\mathfrak{B}}_k =  \big\langle \nabla^2 \bar{\omega}  ; \phi_{k} \big\rangle , \quad \widehat{\mathfrak{L}}_{i,k} = \big\langle \nabla^2 \phi_{i}; \phi_{k} \big\rangle.
\end{align*}}

It remains to compute or assume a good estimate for $\nu_e$. Using an a priori estimate for $\nu_e$ can produce a stable ROM solution. However, as the flow evolves, this prior value might become less effective. Therefore, there should be a strategy to dynamically update this estimate based on the flow conditions/regimes. 

In this regard, we borrow ideas from meteorological data assimilation to correct and update our parameter estimate using live and realistic (possibly noisy) measurements. In particular, we use the forward sensitivity method (FSM), described in Section~\ref{sec:fsm}, to compute an optimal value for eddy viscosity given a few field observations. This also allows us to update our estimate whenever a new observation becomes available. We start with a prior estimate of eddy viscosity (e.g., zero if no priors are available), and solve the ROM equation for a given period of time, $T_w$. As we solve GROM, we also collect some field measurements during this period $T_w$. A penalty term is thus computed as the difference between the GROM prediction and observations, which is used to update our prior estimate for $\nu_e$. This updated value is therefore used to evolve the GROM until new observations become available to match with model's predictions, and so on. The period over which measurements are collected $T_w$ is called the data assimilation window. Note that model's states (e.g., $a_k(t)$) can be different from the measured quantities (e.g., $u(x,t)$), and a mapping between model space and observation space has to be defined. In the following, we formalize our framework for FSM-based eddy viscosity estimation for GROM, called GROM-FSM in the present study. Defining our dynamic model as
\begin{align}
    \dfrac{\mathrm{d} \mathbf{a} }{\mathrm{d} t} = \mathbf{f}(\mathbf{a},\nu_e),
\end{align}
where $\mathbf{a}$ is the vector of modal coefficients defined as $\mathbf{a} = [a_1, a_2, \dots, a_R]^T$ (the superscript $T$ denotes transpose). The (time-continuous) model map $\mathbf{f} = [f_1, f_2, \dots, f_R]^T$ is defined as \textcolor{rev}{$f_k = \mathfrak{B}_k + \nu_e \widehat{\mathfrak{B}}_k +  \sum_{i=1}^{R} \mathfrak{L}_{i,k} a_i + \nu_e \sum_{i=1}^{R} \widehat{\mathfrak{L}}_{i,k} a_i + \sum_{i=1}^{R} \sum_{j=1}^{R} \mathfrak{N}_{i,j,k} a_i a_j$}.

A time-discretization scheme can be utilized to convert this model from continuous-time map $\mathbf{f}$ to a discrete-time map $\mathbf{M}$ as
\begin{equation}
    \mathbf{a}^{k+1} = \mathbf{M}(\mathbf{a}^{k},\nu_e),
\end{equation}
where the superscript $k$ denotes the time index. In our implementation, we adopt the fourth-order Runge-Kutta scheme (RK4) for temporal discretization.

Suppose we collect measurements $\mathbf{z}^{k}$ at a single time instant $t_{k}$, where $t_{k} \in [0,T_w]$. The forecast error is defined at $t_{k}$ as 
\begin{equation}
    \mathbf{e}_{F}^{k} = \mathbf{z}^{k} - \mathbf{h}(\mathbf{a}^{k}),
\end{equation}
where $\mathbf{h}(\cdot)$ defines the mapping from model space to observation space. In our results, we consider two mapping cases. In the first case, we preprocess field observations to compute the ``observed'' coefficients (i.e., $\mathbf{z}^{k} =  \mathbf{a}_{Obs}^{k}$), where the mapping is simply identity matrix (i.e., $\mathbf{h}(\mathbf{a}^{k}) = \mathbf{a}^{k}$). In the second case, we keep observations as velocity field measurement ($\mathbf{z}^{k} =  \mathbf{u}_{Obs}^{k}$), where the mapping becomes a reconstruction map (i.e., $\mathbf{h}(\mathbf{a}^{k}) = \mathbf{u}^{k}$). Specific details are to be given in Section~\ref{sec:results}. 

Although the FSM can be used to treat uncertainties in initial conditions as well as model parameters, we only consider the estimation of the eddy viscosity parameter $\nu_e$. Thus, 
\begin{equation} \label{eq:rom_fsm1}
    \mathbf{H}_2^{k} \delta \nu_e = \mathbf{e}_{F}^{k},
\end{equation}
where $\mathbf{H}_2^{k} = \mathbf{D}^{k}_{\mathbf{a}}(\mathbf{h}) \mathbf{V}^{k}$ as defined in Section~\ref{sec:fsm}. Details of defining model Jacobian are given in \textbf{Appendix A}. For more than a single observation time, we stack Eq.~\ref{eq:rom_fsm1} at different observation times to get the following equation,
\begin{equation} \label{eq:rom_fsm2}
    \mathbf{H}_2 \delta \nu_e = \mathbf{e}_{F}.
\end{equation}
Also, a block-diagonal matrix $\mathbf{R}$ is constituted with the measurement covariance matrices $\mathbf{R}^{k}$ at subsequent observation times. Equation~\ref{eq:rom_fsm2} defines an over-determined system of linear equations in $\delta \nu_e$. A weighted least-squares solution can be computed, with a weighting matrix of $\mathbf{R}^{-1}$ as follows,
\begin{equation}
   \delta \nu_e = \left( \mathbf{H}_2^T \mathbf{R}^{-1} \mathbf{H}_2 \right)^{-1} \mathbf{H}_2^T \mathbf{R}^{-1} \mathbf{e}_{F}, \label{eq:FSMnue}
\end{equation}
where $\delta \nu_e$ is added to our prior estimate of $\nu_e$ (also called background)  to obtain a better approximation and the process is repeated until convergence. The procedure for using FSM to compute the eddy viscosity is summarized in Algorithm~\ref{alg:FSMc}. A tolerance limit has to be set to define convergence (e.g., $1\times 10^{-6}$). We also note that an initial guess for eddy viscosity parameter is required for proper implementation of the algorithm. If no prior knowledge of $\nu_e$ is available, a zero initial guess usually works fine. \textcolor{rev}{Meanwhile, since collected FOM snapshots are already available during an offline stage, they can be treated as field measurement data with negligible noise (corresponding to the underlying solution's assumptions and numerical approximations). Thus, Algorithm~\ref{alg:FSMc} can be applied offline along with the construction of GROM model to provide a prior estimate of the suitable closure parameter.} 

\begin{algorithm}[ht]
\SetAlgoLined
\SetKwInOut{Input}{Input}\SetKwInOut{Output}{Output}
\Input{Dynamic model $\mathbf{M}(\cdot)$, observation operator $\mathbf{h}(\cdot)$, initial condition $\mathbf{a}^1$, a set of observations $\mathbf{z}^1, \mathbf{z}^2, \dots \mathbf{z}^N$, an initial guess for eddy viscosity parameter $\nu_e$, and a tolerance $tol$ value}
\Output{An estimate of the eddy viscosity $\nu_e$}
\medskip
initialization \\
\For{$ i \leftarrow 1$ \KwTo $max \ iter$}{
    $\mathbf{V}^1 = \mathbf{0}$ \\
    $\mathbf{e}_{F} = \mathbf{z}^{1} - \mathbf{h}(\mathbf{a}^{1})$ \\  
    $\mathbf{H}_2 = \mathbf{D}^{1}_{\mathbf{a}}(\mathbf{h}) \mathbf{V}^{1}$ \\
    $\mathbf{R} = \mathbf{R}^{1}$ \\
    \For{$n\leftarrow 1$ \KwTo $N-1$}{
        $\mathbf{a}^{n+1} = \mathbf{M}(\mathbf{a}^{n}, \nu_e)$ \\
        $\mathbf{V}^{n+1} = \mathbf{D}_{\mathbf{a}}^n
        (\mathbf{M})\mathbf{V}^{n} + \mathbf{D_{\nu_e}^n (\mathbf{M})} $ \\
        \If{(observation $\mathbf{z}^{n+1}$ is available)}{
               $\mathbf{e}_{F}^{n+1} = \mathbf{z}^{n+1} - \mathbf{h}(\mathbf{a}^{n+1})$ \\                
               $\mathbf{H}_2^{n+1} = \mathbf{D}^{n+1}_{\mathbf{a}}(\mathbf{h}) \mathbf{V}^{n+1}$ \\
               $\mathbf{e}_{F} = \begin{bmatrix} 
                \mathbf{e}_{F} \\ \mathbf{e}_{F}^{n+1}  \end{bmatrix}$, \qquad
               $\mathbf{H}_2 = \begin{bmatrix}
                \mathbf{H}_2 \\ \mathbf{H}_2^{n+1}  \end{bmatrix}$ \\
               $\mathbf{R} = \begin{bmatrix}
                \mathbf{R} & \\ & \mathbf{R}^{n+1}  \end{bmatrix}$
        }       

    }
    $\delta \nu_e = \left( \mathbf{H}_2^T \mathbf{R}^{-1} \mathbf{H}_2 \right)^{-1} \mathbf{H}_2^T \mathbf{R}^{-1} \mathbf{e}_{F}$  \\
    \eIf{($\delta \nu_e \le tol$)}{
        $break$ \\
    }{
        $\nu_e = \nu_e + \delta \nu_e$ \\
    }
}
 \caption{Forward sensitivity method for estimating eddy viscosity in GROM closure }
 \label{alg:FSMc}
\end{algorithm}


\section{Results} \label{sec:results}
\begin{figure}[ht]
\centering
\includegraphics[width=0.95\linewidth]{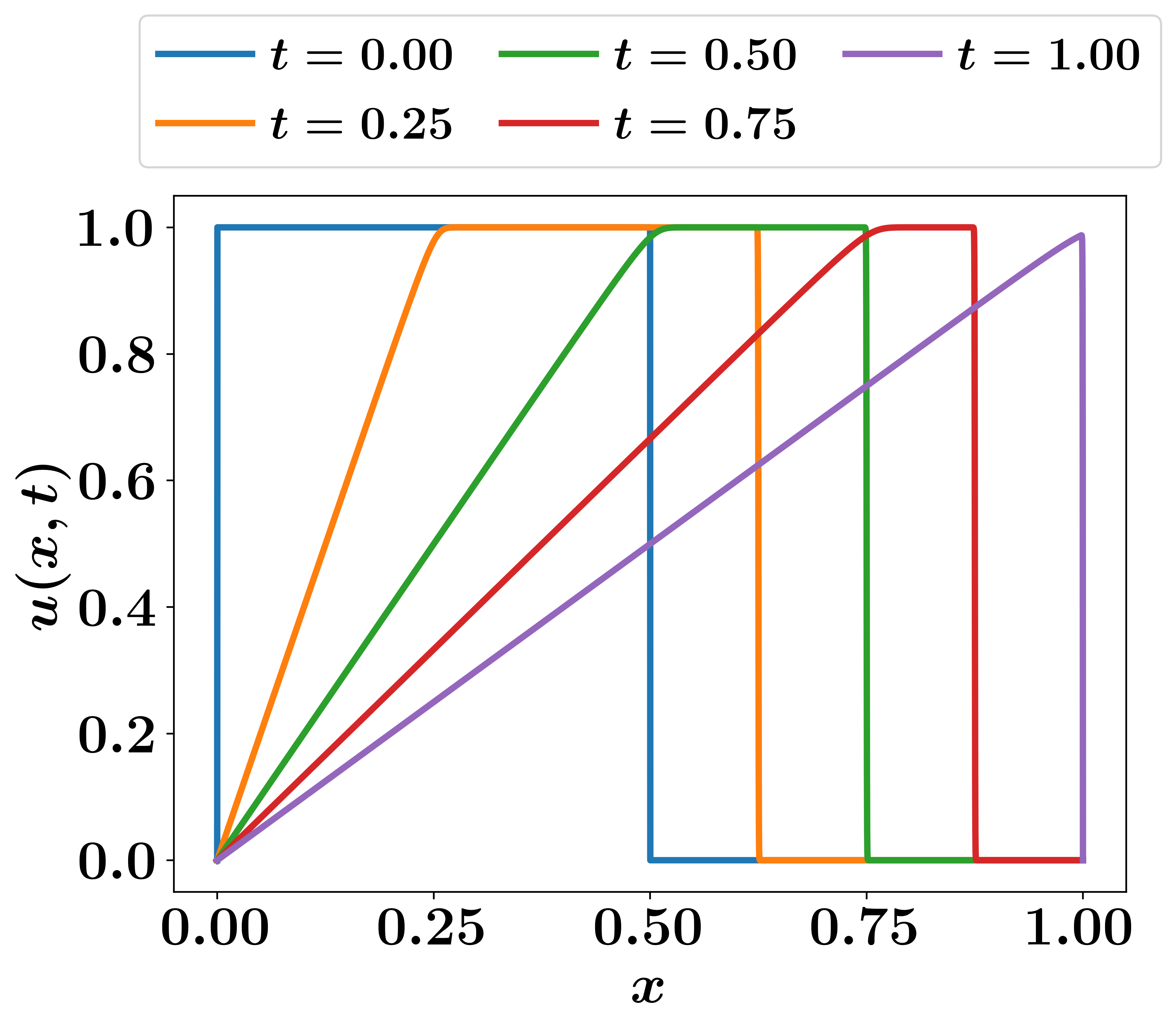}
\caption{Evolution of the FOM velocity field, characterized by a moving shock with square wave.}
\label{fig:FOM}
\end{figure}

\textcolor{rev}{In this section, we present our results for the utilization the proposed methodology to compute and update the eddy viscosity parameter via FSM, applied to the introduced two test problems (i.e., 1D Burgers problem and 2D Kraichnan turbulence).
\subsection{1D Burgers Problem}}
For the 1D Burgers problem, we assume an initial condition of a square wave defined as
\begin{equation}
u(x,0) = \begin{cases}
            1 , &\quad\text{if } 0 < x \le L/2 \\
            0 , &\quad\text{if } L/2 < x \le L,
        \end{cases}
\end{equation}
with zero Dirichlet boundary conditions, $u(0,t) = u(L,t) = 0$. We consider a spatial domain of $L=1$, and solve at $\text{Re} = 10^4$ for $t \in [0,1]$. For numerical computations, we use a family of fourth order compact schemes for spatial derivatives \cite{lele1992compact}, and skew-symmetric formulation for the nonlinear term. Also, we use the fourth order Runge-Kutta (RK4) scheme for temporal integration with a time step of $10^{-4}$ over a spatial grid of $4096$. For POD basis generation, we collect 100 snapshots (i.e., every 100 time steps). The temporal evolution of the 1D Burgers problem using the described setup is shown in Figure~\ref{fig:FOM}, where we can see the advection of the shock wave.

The described Burgers problem with square wave is challenging for ROM applications. In our GROM implementation, we consider $R=8$ modes and $\Delta t = 0.01$ for time integration. In the following, we discuss the estimation of eddy viscosity via FSM using full and sparse field measurements.

\subsubsection{Full field measurement} \label{sec:full}

\begin{figure}[ht]
\centering
\includegraphics[width=0.95\linewidth]{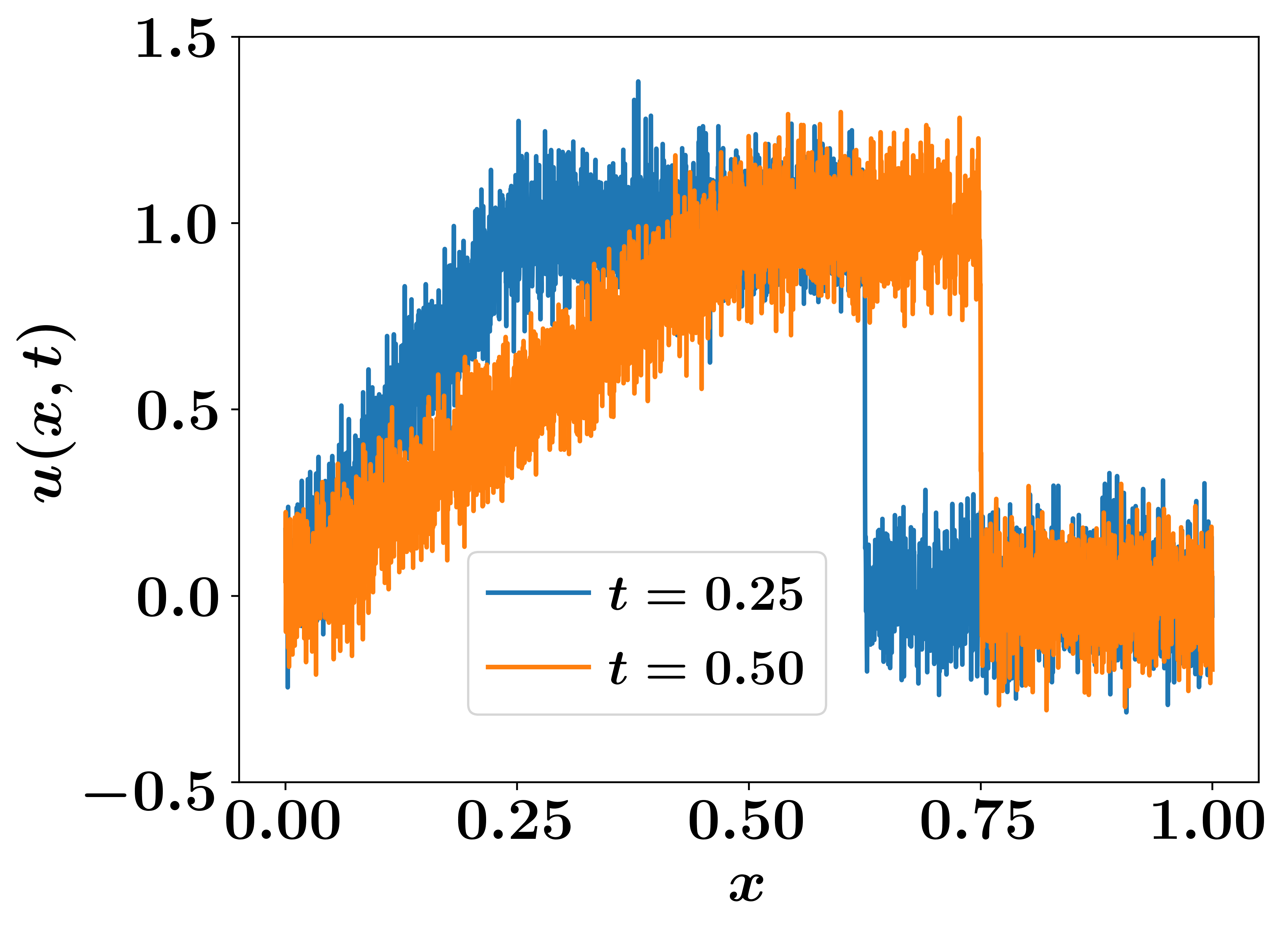}
\caption{Noisy measurement of velocity fields at $t=0.25$ s and $t=0.50$ s, assuming sensors are located at all grid points.}
\label{fig:FullObs}
\end{figure}

In our first case, we investigate the assimilation of noisy full field measurement as
\begin{equation}
    u_{Obs} (x,t) = u(x,t)  + v(x,t),
\end{equation}
where $v(x,t)$ is a white Gaussian noise with zero mean and covariance matrix $\mathbf{R}(t)$. In particular, we define $\mathbf{R}(t) = \sigma_{Obs}^2 \mathbf{I}$, with $\sigma_{Obs} = 0.1$. We assume a data assimilation window of $0.5$ s and collect measurements at $t=0.25$ and $t=0.5$, as demonstrated in Figure~\ref{fig:FullObs}.

\begin{figure}[ht]
\centering
\includegraphics[width= 0.95\linewidth]{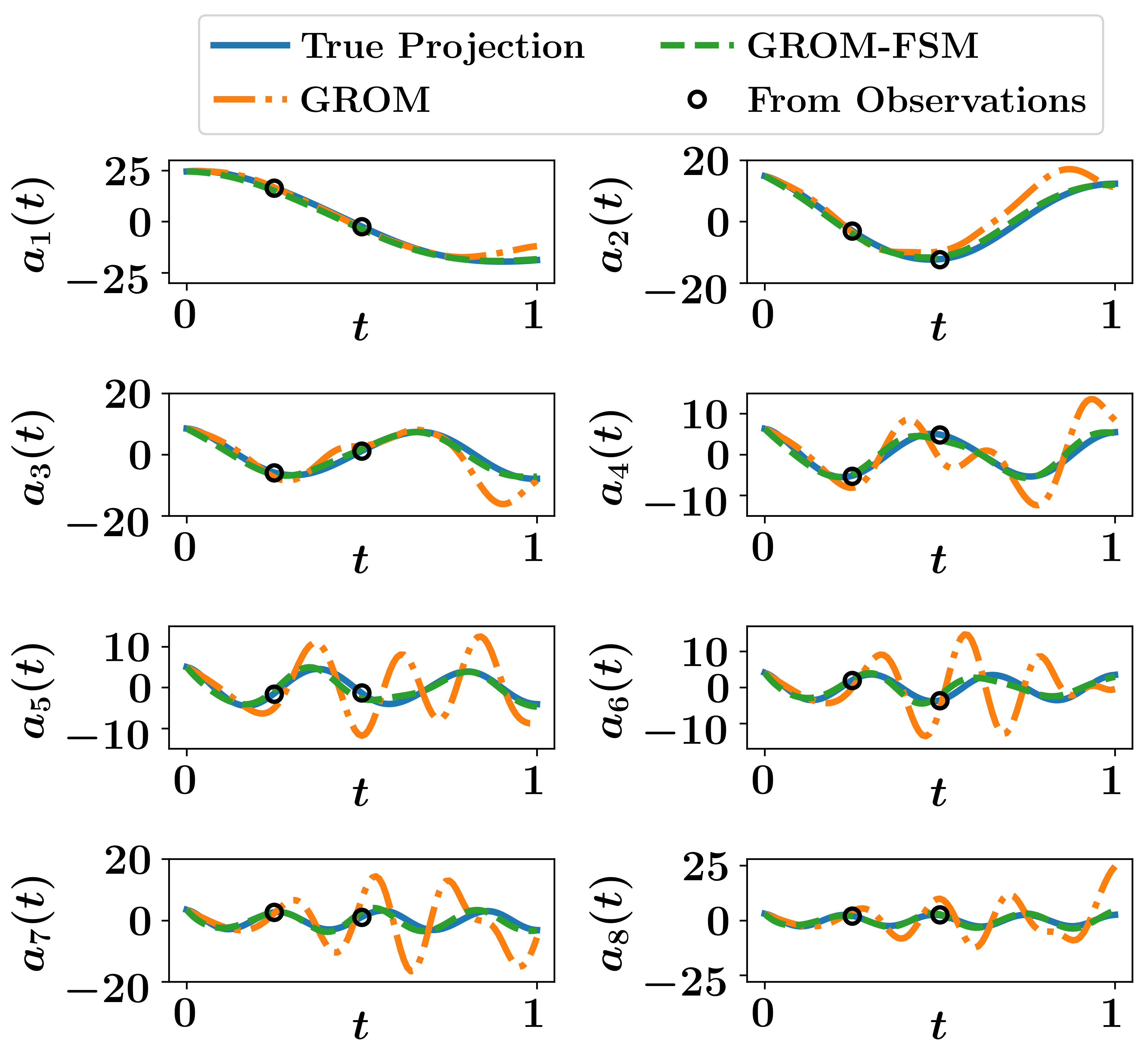}
\caption{Temporal evolution of POD coefficients, assuming full field measurements are available.}
\label{fig:aFullObs}
\end{figure}

Instead of defining a map between model space and observation space, we preprocess our measurement by projecting them onto the POD basis to compute the ``observed'' modal coefficients as
\begin{equation}
    a^k_{i,Obs} = \langle \mathbf{u}^k_{Obs} - \bar{\mathbf{u}}; \phi_i \rangle.
\end{equation}
Thus, \textcolor{rev}{$\mathbf{z}^k = \mathbf{a}^k_{Obs}$} and the observation operator is defined $\mathbf{h}(\mathbf{a}) = \mathbf{a}$, with a Jacobian \textcolor{rev}{equal to the identity matrix (i.e.,  $\mathbf{D}_{\mathbf{a}}(\mathbf{h}) = \mathbf{I}_{R}$, where $\mathbf{I}_{R}$ is the $R \times R$ identity matrix)}. Also, the observational covariance matrix is set as $\mathbf{R}^k = \sigma_{Obs}^2 \mathbf{I}_{R}$. If we implement the procedure described in Section~\ref{sec:fsm} to obtain an estimate for $\nu_e$ and solve GROM with and without closure, we obtain the results in Figure~\ref{fig:aFullObs} for the temporal evolution the modal coefficients. For comparison, we also plot the true projection values of $\mathbf{a}$, defined as 
\begin{equation}
    a^k_{i} = \langle \mathbf{u}^k_{FOM} - \bar{\mathbf{u}}; \phi_i \rangle.
\end{equation}

\begin{figure}[ht]
\centering
\includegraphics[width=0.95\linewidth]{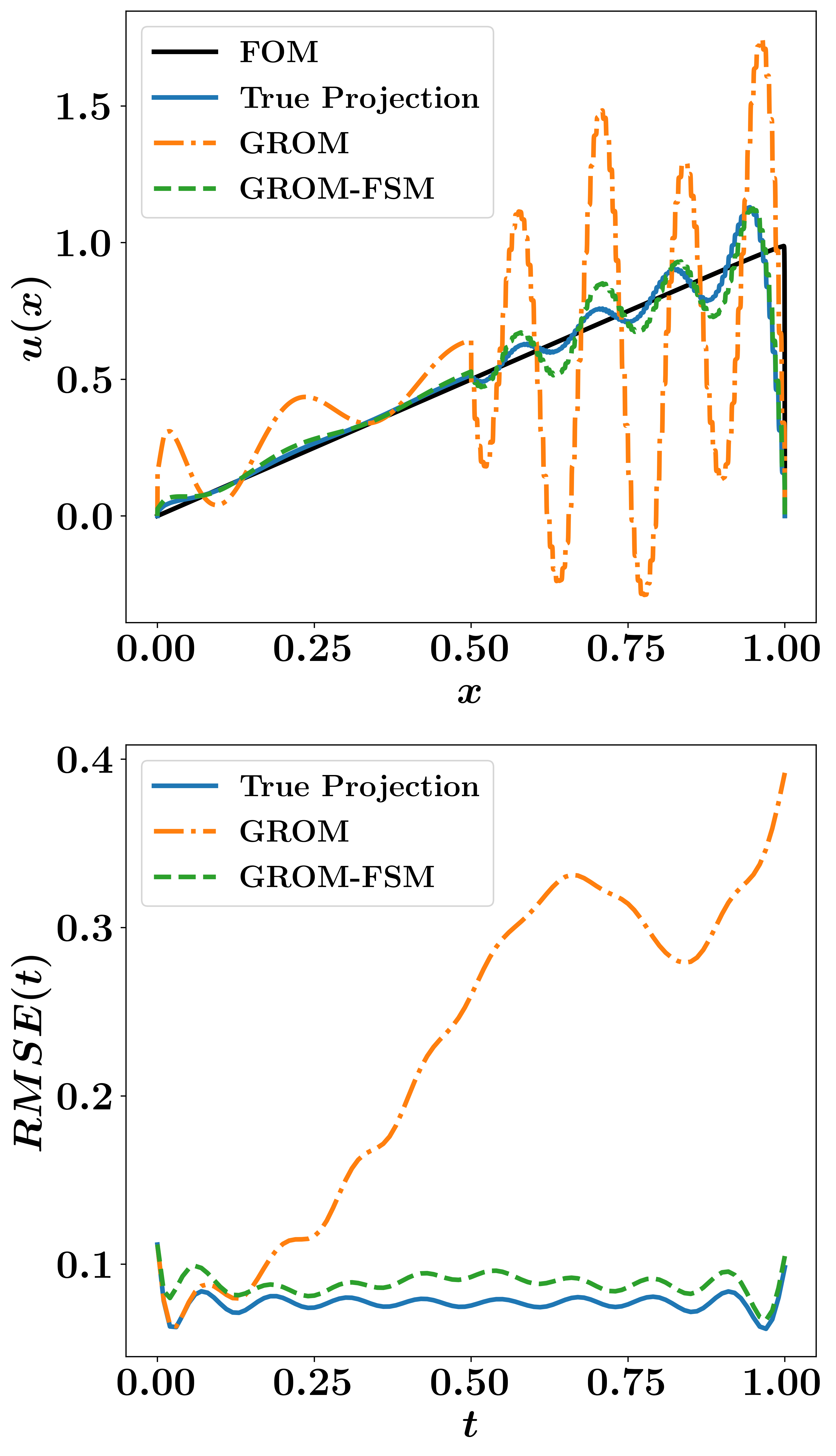}
\caption{Velocity field reconstruction in case of full field measurements. Top: reconstruction of final velocity field using GROM and GROM with FSM eddy viscosity compared to the FOM and true projection fields. Bottom: RMSE of reconstructed fields at different time instants.}
\label{fig:uFullObs}
\end{figure}

Also, we sketch reconstructed velocity field at final time $t=1$ in \textcolor{rev}{Figure~\ref{fig:uFullObs}}. It is clear that GROM without closure is unable to capture the true dynamical behavior of the described Burgers problem. On the other hand, GROM-FSM is shown to almost match the true projection. It is assumed that true projected values represent the best values that projection-based ROM can provide. For quantitative assessment, we also plot the root mean squares error (RMSE) of ROM predictions defined as
\begin{equation}
    RMSE(t) = \sqrt{\dfrac{1}{n}\sum_{i=1}^{n}{\bigg(u_{FOM}(x_i,t) - u_{ROM}(x_i,t) \bigg)^2} }
\end{equation}

\subsubsection{Sparse field measurement}
\begin{figure}[ht]
\centering
\includegraphics[width=0.95\linewidth]{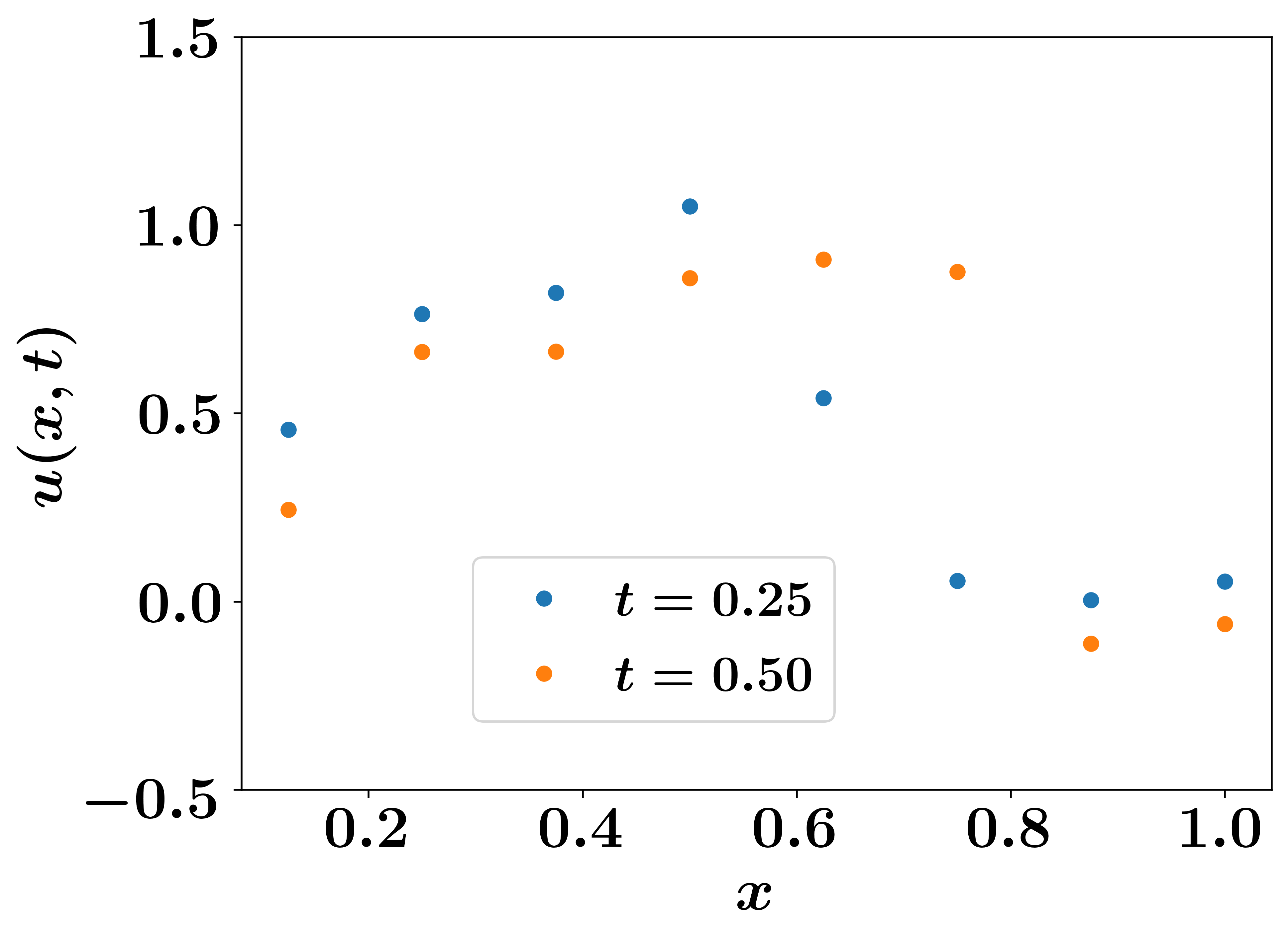}
\caption{Noisy measurement of velocity fields at $t=0.25$ s and $t=0.50$ s, assuming sensors are located at 8 grid points.}
\label{fig:PartObs}
\end{figure}

Since full field measurements are usually inaccessible, we extend our study to consider sparse field measurements. In particular, we locate sensors at 8 points, equally spaced at $1/8, 2/8, 3/8, 4/8, 5/8, 6/8, 7/8, 8/8$ as shown in Figure~\ref{fig:PartObs}. To assimilate those measurements, we consider two cases. The first one is similar to the full field measurement case, where we preprocess those measurements to compute a least-squares approximation of the corresponding observed modal coefficients. In the second case, we keep our observation as field measurements and define an operator to map model state (i.e., modal coefficients) to observations (i.e., velocity).

\begin{figure}[ht]
\centering
\includegraphics[width=0.95\linewidth]{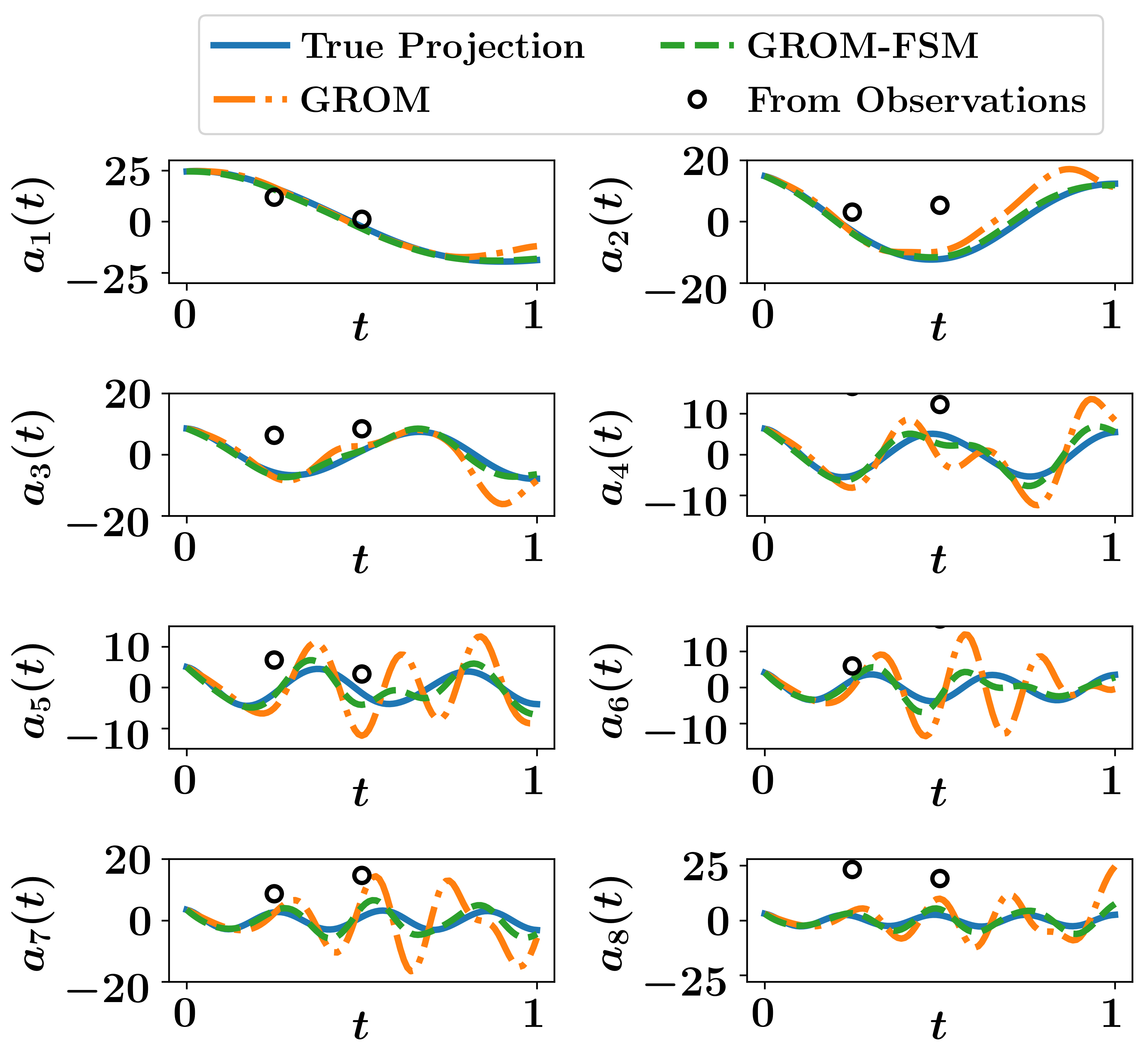}
\caption{Temporal evolution of POD coefficients, where sparse field measurements are preprocessed to estimate the observed POD coefficients.}
\label{fig:aPartObsInv}
\end{figure}

\paragraph{\textbf{From measurements to POD coefficients}} \label{sec:sp1}
In order to preprocess the sparse measurements to approximate the observed modal coefficients, we sample Eq.~\ref{eq:qROM} at the sensors locations as follows,
\textcolor{rev}{
\begin{align}\label{eq:partObs1}
    &\begin{bmatrix}
        \phi_1(x_{O1}) & \phi_2(x_{O1}) & \dots & \phi_R(x_{O1}) \\
        \phi_1(x_{O2}) & \phi_2(x_{O2}) & \dots & \phi_R(x_{O2}) \\
        \vdots         &                &       & \vdots         \\
        \phi_1(x_{O8}) & \phi_2(x_{O8}) & \dots & \phi_R(x_{O8})
    \end{bmatrix}
    \begin{bmatrix} a^k_{1,Obs} \\ a^k_{2,Obs} \\  \vdots \\ a^k_{R,Obs} \end{bmatrix}  \nonumber \\
    = &\begin{bmatrix} u^k_{Obs}(x_{O1}) - \bar{u}(x_{O1}) \\ u^k_{Obs}(x_{O2}) - \bar{u}(x_{O2}) \\  \vdots \\ u^k_{Obs}(x_{O8}) - \bar{u}(x_{O8}) \end{bmatrix},
\end{align}}
which can be generally solved using the pseudo-inverse. Then, the same observation operator and its Jacobian as defined in Section~\ref{sec:full} are used. The temporal evolution of the modal coefficients are given in Figure~\ref{fig:aPartObsInv}. Although the GROM-FSM results are better than GROM, they are significantly worse than those in Figure~\ref{fig:aFullObs}. Of course, this is to be expected since we are using measurements at only 8 points, rather than 4096 locations. However, we also find that the observed modal coefficients calculations using Eq.~\ref{eq:partObs1} is greatly sensitive to the level of noise. Indeed, we find that least-squares computations sometimes do not converge (a remedy will be provided in Section~\ref{sec:sp2}). Moreover, we can see that the POD modal coefficients from observations are significantly different than the true ones.

\begin{figure}[ht]
\centering
\includegraphics[width=0.95\linewidth]{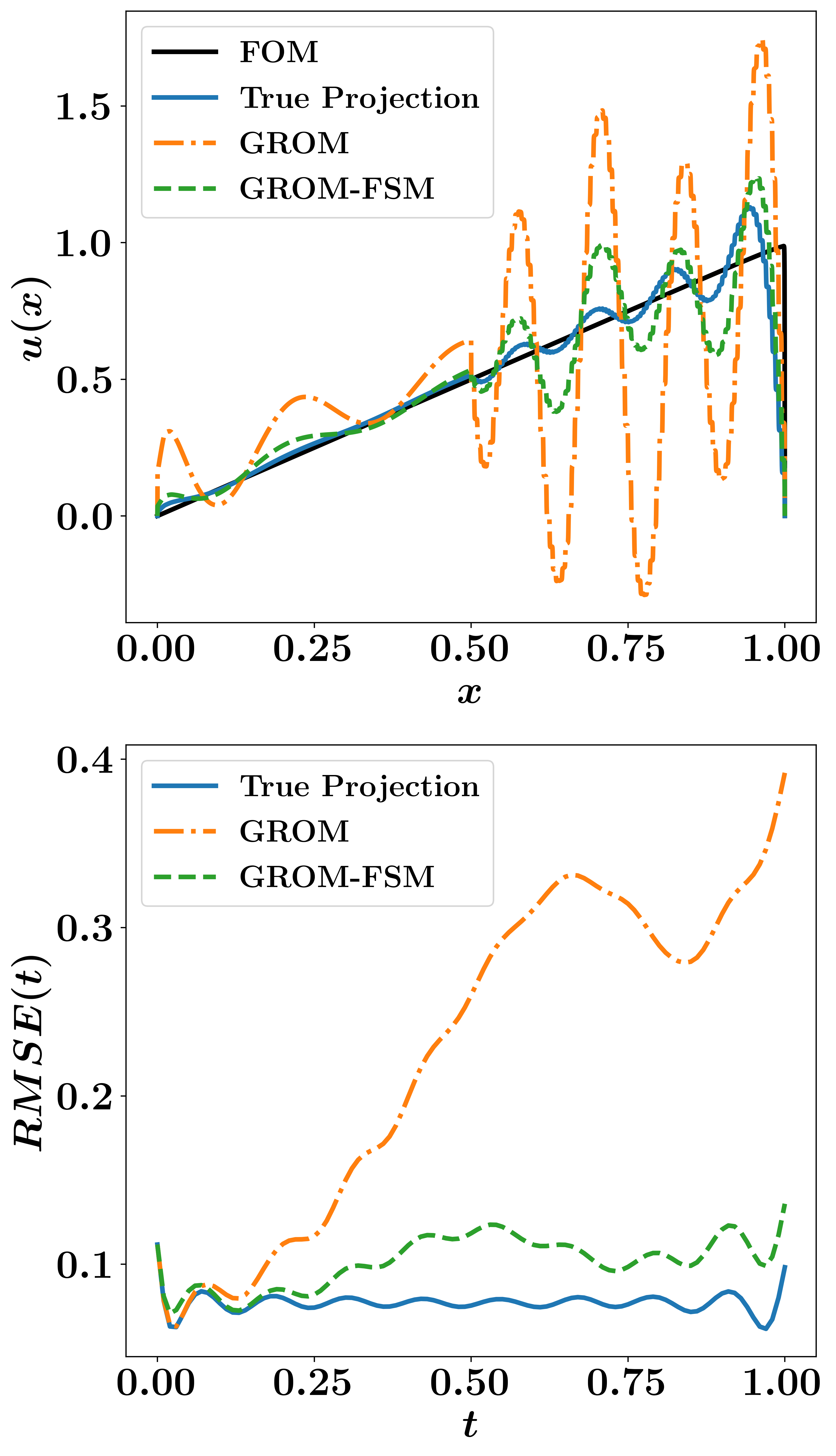}
\caption{Velocity field reconstruction in case of preprocessing sparse field measurements to compute the observed POD coefficients. Top: reconstruction of final velocity field using GROM and GROM with FSM eddy viscosity compared to the FOM and true projection fields. Bottom: RMSE of reconstructed fields at different time instants.}
\label{fig:uPartObsInv}
\end{figure}

The reconstructed field at final time as well as the $RMSE$ at different times are demonstrated in Figure~\ref{fig:uPartObsInv}. We see that a small improvement is obtained in GROM-FSM, compared to GROM. We also note that for different noise levels, we get different performances for the GROM-FSM. This implies that this way of assimilating sparse observations is less reliable, and a more robust approach should be utilized. In Section~\ref{sec:sp2}, we discuss another way of using sparse observations to perform data assimilation for ROM closure.

\paragraph{\textbf{From POD coefficients to measurements}} \label{sec:sp2}
Now, we discuss defining an observational operator to construct a robust map between model space and measurement space. Similar to Section~\ref{sec:sp1}, we sample Eq.~\ref{eq:qROM} at sensor location, but we introduce a map to reconstruct the velocity field at these locations using the model predicted coefficients. In other words, in Section~\ref{sec:sp1}, we use the sensors measurements to approximate a value for $\mathbf{a}^k_{Obs}$. But in this section, we use model predicted coefficients $\mathbf{a}^k$ to approximate the velocity field values at sensor locations (i.e., $ u^k(x_{O1}), u^k(x_{O2}), \dots, u^k(x_{O8})$) as follows,
\begin{align}\label{eq:partObs2}
    &\begin{bmatrix}
        \phi_1(x_{O1}) & \phi_2(x_{O1}) & \dots & \phi_R(x_{O1}) \\
        \phi_1(x_{O2}) & \phi_2(x_{O2}) & \dots & \phi_R(x_{O2}) \\
        \vdots         &                &       & \vdots         \\
        \phi_1(x_{O8}) & \phi_2(x_{O8}) & \dots & \phi_R(x_{O8})
    \end{bmatrix}
    \begin{bmatrix} a^k_{1} \\ a^k_{2} \\  \vdots \\ a^k_{R} \end{bmatrix} \\
    =  &\begin{bmatrix} u^k(x_{O1}) - \bar{u}(x_{O1})\\ u^k(x_{O2}) - \bar{u}(x_{O2}) \\  \vdots \\ u^k(x_{O8}) - \bar{u}(x_{O8}) \end{bmatrix}.
\end{align}

\begin{figure}[ht]
\centering
\includegraphics[width=0.95\linewidth]{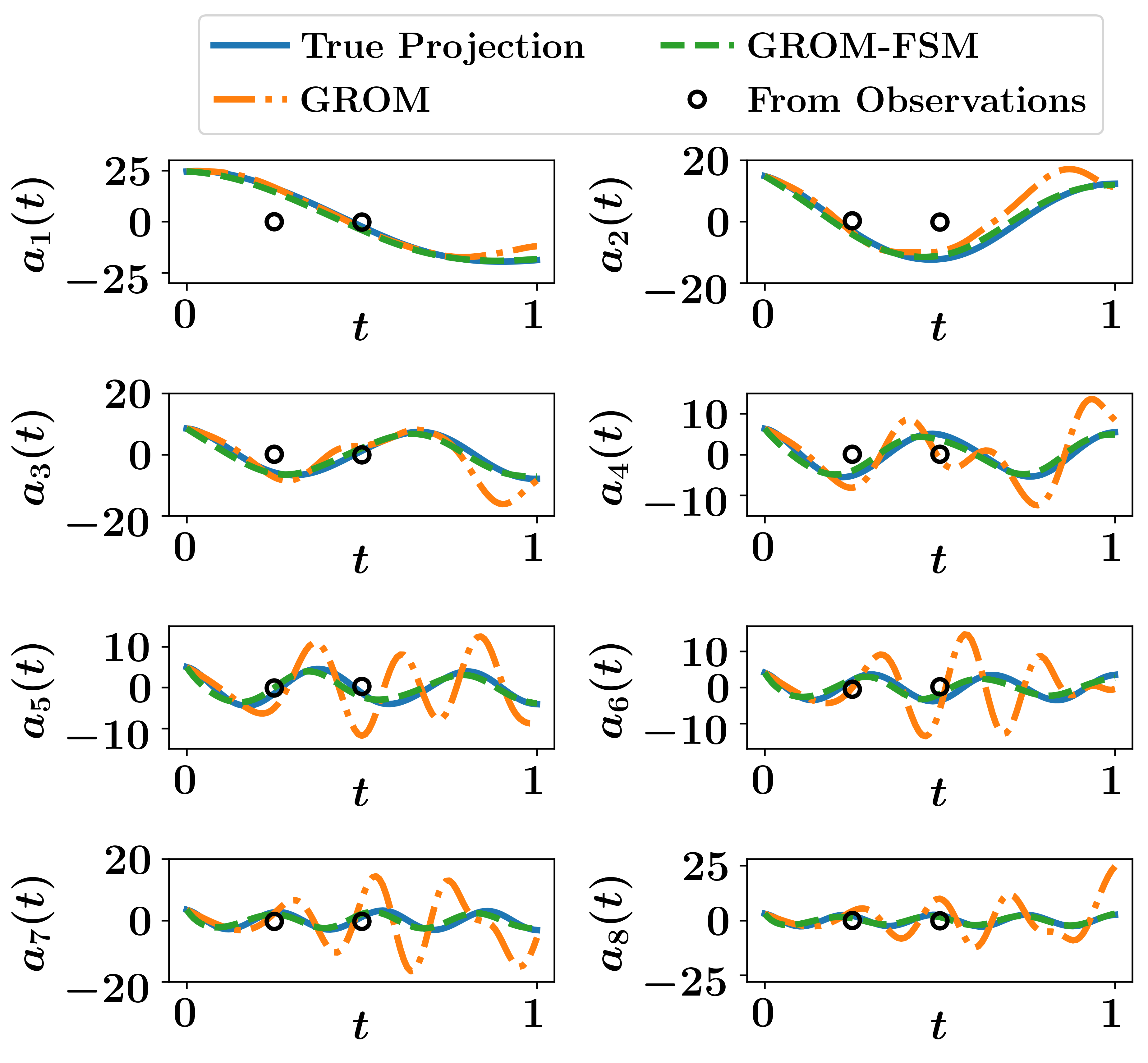}
\caption{Temporal evolution of POD coefficients, where sparse field measurements are compared against POD field reconstruction using the observer operator $\mathbf{C}$.}
\label{fig:aPartObsRec}
\end{figure}

\begin{figure}[ht]
\centering
\includegraphics[width=0.95\linewidth]{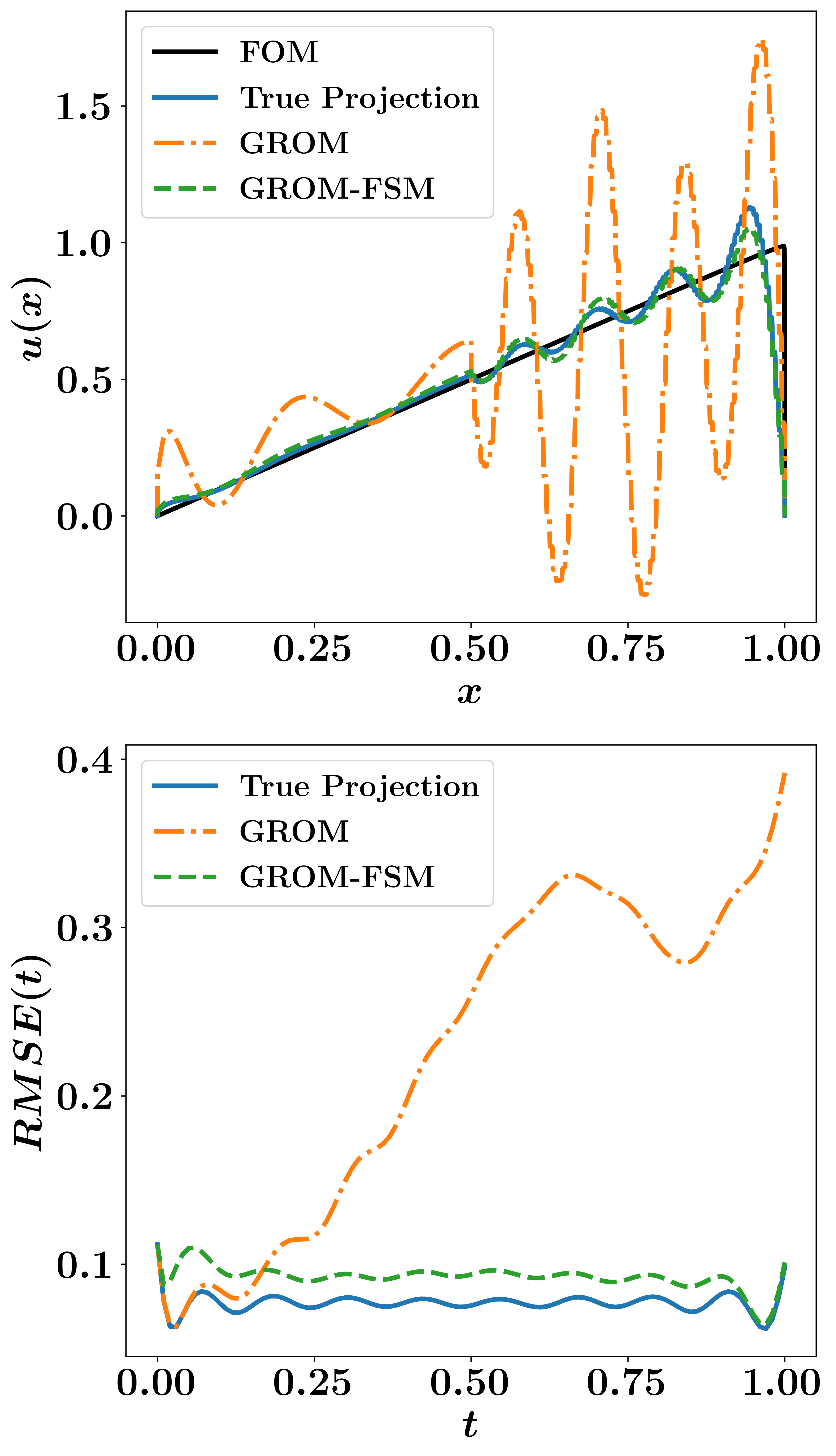}
\caption{Velocity field reconstruction where sparse field measurements are compared against POD field reconstruction using the observer operator $\mathbf{C}$. Top: reconstruction of final velocity field using GROM and GROM with FSM eddy viscosity compared to the FOM and true projection fields. Bottom: RMSE of reconstructed fields at different time instants.}
\label{fig:uPartObsRec}
\end{figure}

Thus, we define $\mathbf{z}^k = \mathbf{u}^k_{Obs}$, and the observation operator $\mathbf{h}(\mathbf{a}) = \mathbf{C} \mathbf{a}$, where $\mathbf{C}$ is the matrix of basis functions sampled at sensors locations as follows,
\begin{equation}\label{eq:mapsp2}
    \mathbf{C} = \begin{bmatrix}
        \phi_1(x_{O1}) & \phi_2(x_{O1}) & \dots & \phi_R(x_{O1}) \\
        \phi_1(x_{O2}) & \phi_2(x_{O2}) & \dots & \phi_R(x_{O2}) \\
        \vdots         &                &       & \vdots         \\
        \phi_1(x_{O8}) & \phi_2(x_{O8}) & \dots & \phi_R(x_{O8})
    \end{bmatrix}.
\end{equation}
Thus, the Jacobian of $\mathbf{h(\cdot)}$ is defined as $\mathbf{D}_{\mathbf{a}}(\mathbf{h}) = \mathbf{C}$. We repeat the same GROM-FSM implementation with those redefined operators. Results are shown in Figure~\ref{fig:aPartObsRec} and Figure~\ref{fig:uPartObsRec}, where we can see that this approach of assimilating measurements is more robust than the one discussed in Section~\ref{sec:sp1} with higher accuracy. We also note that similar performance is achieved using higher level of noise in measurements, while the approach in Section~\ref{sec:sp1} requires very low level of observational noise.

\begin{figure*}[ht]
\centering
\includegraphics[trim= 60 0 60 0, clip, width=0.9\textwidth]{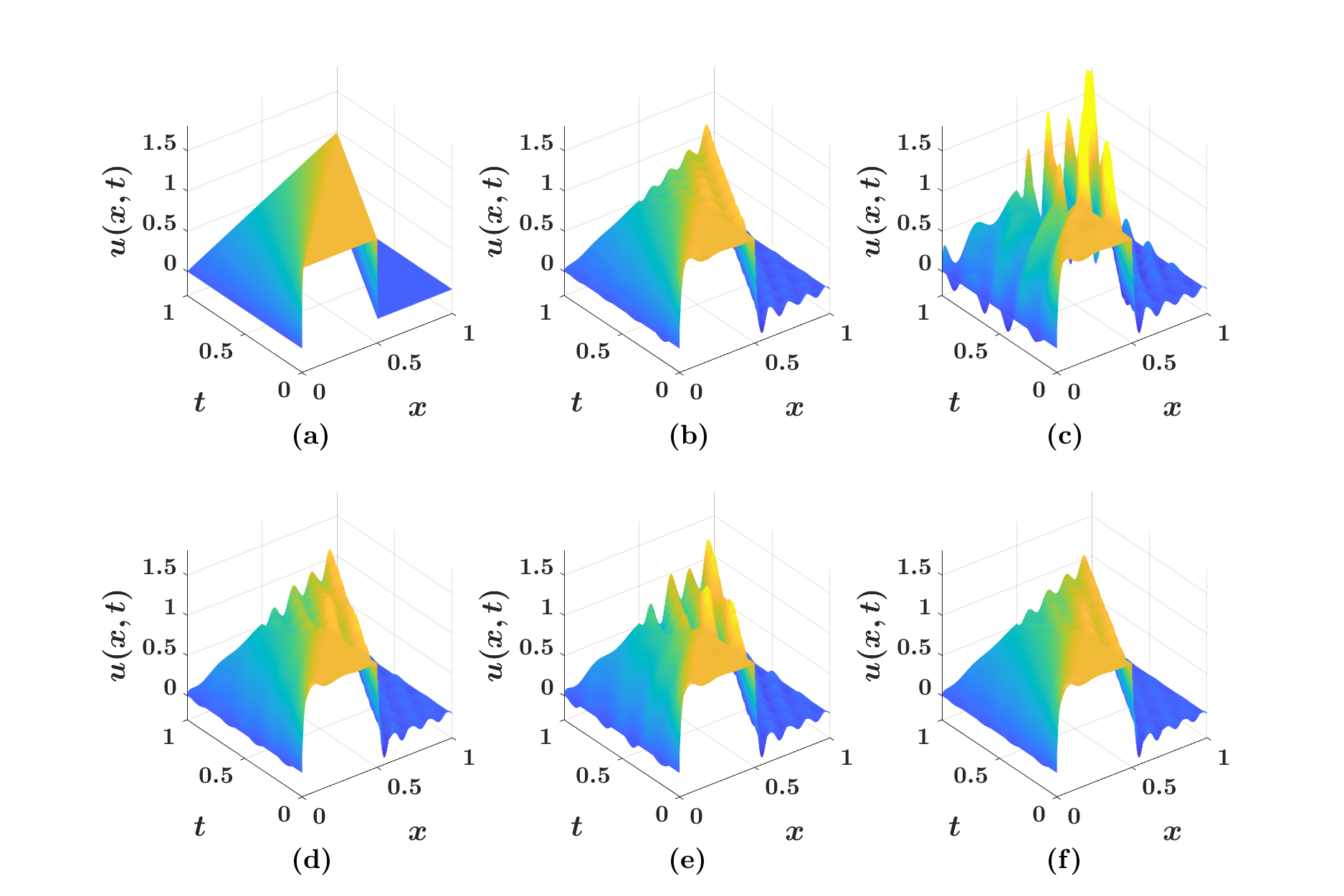}
\caption{\textcolor{rev}{Surface plots for the temporal evolution of velocity fields from (a) FOM, (b) true projection, (c) GROM and FSM with (d) full and (e--f) sparse measurements configurations. Note: (e) shows the results using the method presented in Section~\ref{sec:sp1}, while (f) refers to the method presented in Section~\ref{sec:sp2}.}}
\label{fig:surf}
\end{figure*}

Finally, for a big picture comparison, we plot the spatio-temporal evolution of reconstructed velocity fields for all discussed measurement setups compared to FOM and true projection fields in Figure~\ref{fig:surf}. From this figure, we notice that solution of GROM without closure is unstable, and brings non-physical predictions. On the other hand, predictions of GROM-FSM with full field measurements almost match the true projected fields. Also, assimilating sparse observations via the reconstruction map $\mathbf{C}$ is significantly superior to approximating observed coefficients using the pseudo-inverse approach. The latter shows some non-physical predictions, similar to GROM.

\begin{figure*}[ht]
\centering
\includegraphics[width=0.95\textwidth]{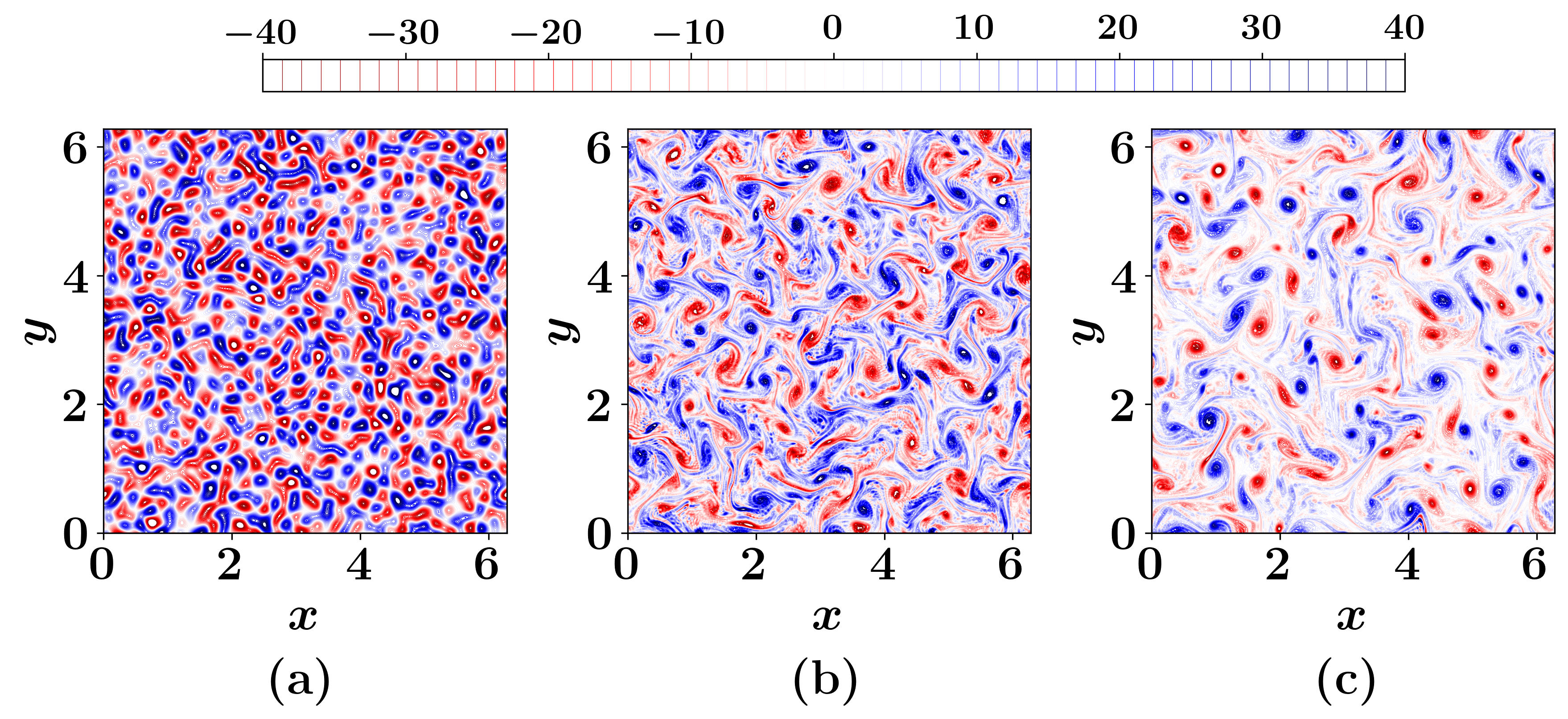}
\caption{\textcolor{rev}{FOM results showing the time evolution of vorticity fields for the 2D turbulence problem: (a) $t=0.0$, (b) $t=2.0$, and (c) $t=4.0$.}}
\label{fig:2DFOM}
\end{figure*}

\textcolor{rev}{\subsection{2D Kraichnan Turbulence}
For 2D turbulence, the inertial range in the energy spectrum is proportional to $k^{-3}$ in the inviscid limit according to the Kraichnan–Batchelor–Leith (KBL) theory \cite{kraichnan1967inertial,batchelor1969computation,leith1971atmospheric}. In our numerical experiments, the initial energy spectrum in Fourier space is given by
\begin{equation}
    E(k) = \dfrac{4k^4}{3\sqrt{\pi} k_p^5}  \exp\left[-\left(\dfrac{k}{k_p}\right)^2\right],
\end{equation}
where $k=\sqrt{k_x^2+k_y^2}$  and $k_p$ is the wavenumber at which the maximum value of initial energy spectrum occurs. During the time evolution process, due to the nonlinear interactions, this spectrum quickly approaches toward $k^{-3}$ spectrum. The magnitude of the vorticity Fourier coefficients is related to the energy spectrum as
\begin{equation}
    |\tilde{\omega}(k)| = \sqrt{\dfrac{k}{\pi} E(k)}.
\end{equation}
Thus, the initial vorticity distribution (in Fourier space) is obtained by introducing a random phase. For more details regarding derivation of the initial vorticity distribution from an assumed energy spectrum can be found in \cite{san2012high}. In the present study, we use a spatial computational domain of $(x,y) \in [0,2\pi] \times [0,2\pi]$ and a time domain of $t\in [0,4]$. Periodic boundary conditions are applied in both $x$ and $y$ directions. A spatial grid of $512^2$ and a timestep of $\Delta t=0.001$ are used for FOM solution, and $800$ snapshots of vorticity fields are stored for POD basis generation. A fourth-order accurate Arakawa scheme \cite{arakawa1966computational} is adopted for spatial discretization. Contours of voriticy field at different time instants are shown in Figure~\ref{fig:2DFOM} initiated by introducing a random phase shift in the Fourier space, with $k_p=10$.}

\textcolor{rev}{For ROM implementation, we consider $R=16$ corresponding to a RIC value of around $80\%$. Similar to the 1D Burgers problem, we test the FSM capabilities to estimate an optimal value of eddy viscosity considering full and sparse field measurements with $\sigma_{Obs}=0.1$. We assume a data assimilation window of $2$, and measurement data are collected at $t=1$ and $t=2$, while testing is performed up to $t=4$. We highlight here that all the 2D fields are rearranged into 1D column vectors to follow the same notations provided in Section~\ref{sec:rom} (e.g., the Euclidean inner product). However, for contour plots, they are reshaped back into 2D fields.}

\begin{figure}[ht]
\centering
\includegraphics[width=0.9\linewidth]{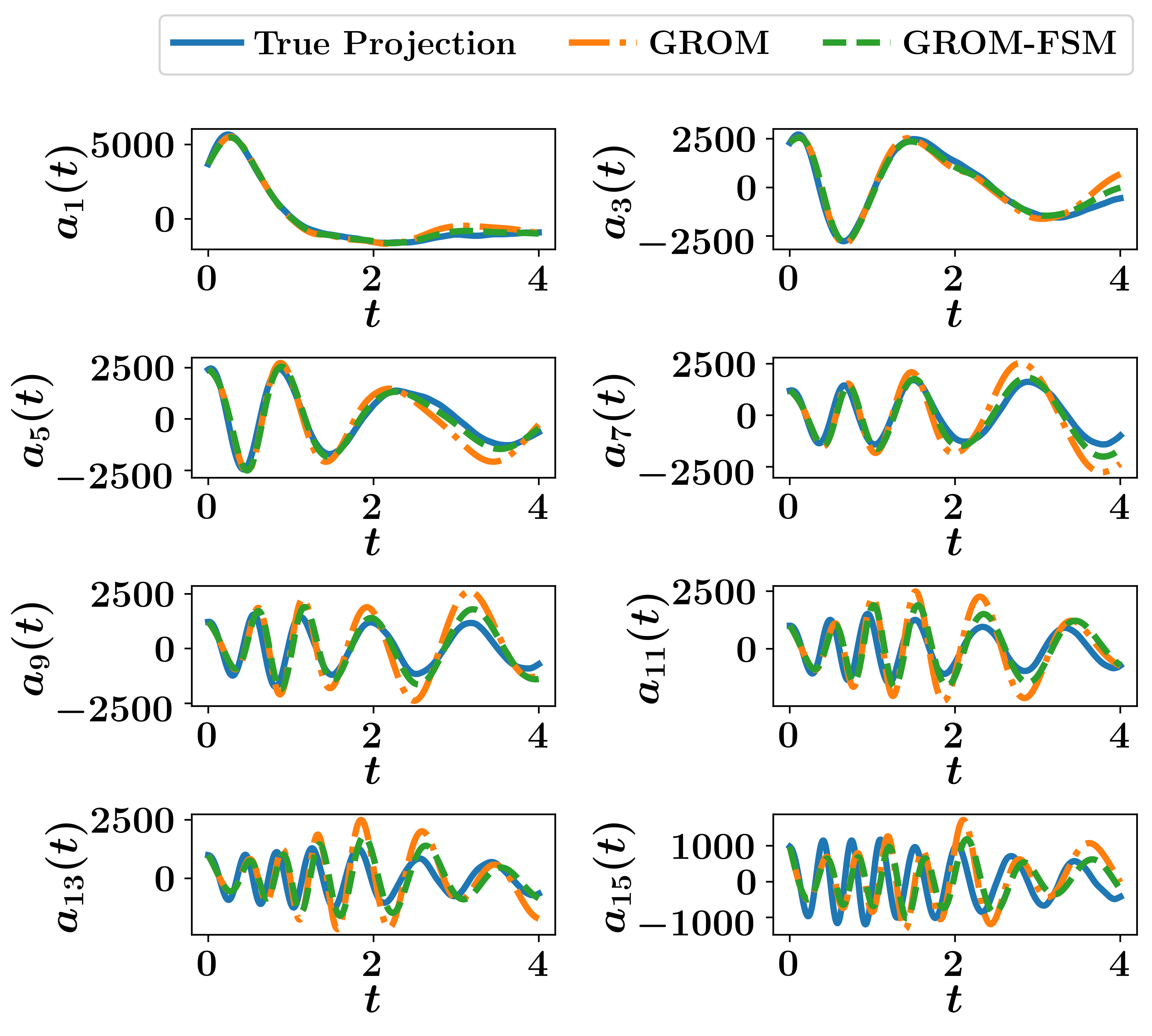}
\caption{\textcolor{rev}{Time evolution of the modal coefficients for the 2D turbulence case, assuming full field measurements are available at $t=1$ and $t=2$.}}
\label{fig:2DaF}
\end{figure}

\textcolor{rev}{\subsubsection{Full field measurement}}
\textcolor{rev}{Since full field measurements are available (though noisy), we can project these field data onto the basis functions $\phi$ to obtain the corresponding \emph{observed} modal coefficients as
\begin{equation}
    a^k_{i,Obs} = \langle \boldsymbol{\omega}^k_{Obs} - \bar{\boldsymbol{\omega}}; \phi_i \rangle.
\end{equation}}
\textcolor{rev}{Following the same procedure as in Section~\ref{sec:full}, an optimal value of eddy viscosity parameter is estimated with the FSM methodology. In Figure~\ref{fig:2DaF}, we show the GROM solution equipped by an eddy viscosity closure, computed by the proposed approach compared to the background solution of standard GROM without closure. Also, we plot the true projection (TP) results, where the true POD modal coefficients are obtained as 
\begin{equation}
    a^k_{i} = \langle \boldsymbol{\omega}^k_{FOM} - \bar{\boldsymbol{\omega}}; \phi_i \rangle.
\end{equation}}

\begin{figure*}[ht]
\centering
\includegraphics[width=0.95\linewidth]{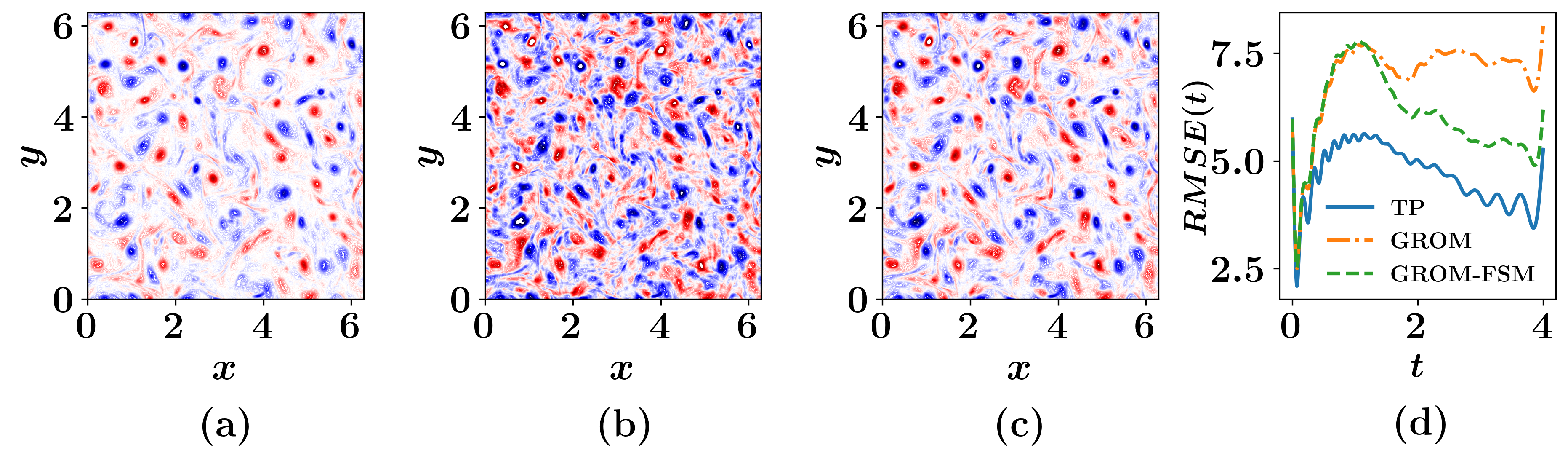}
\caption{\textcolor{rev}{Reconstructed vorticity fields for the 2D turbulence problem at $t=4$ for (a) true projection, (b) GROM, and (c) GROM-FSM along with the $RMSE$ variation with time (d), assuming full field measurements are available at $t=1$ and $t=2$.}}
\label{fig:2DwF}
\end{figure*}

\textcolor{rev}{From Figure~\ref{fig:2DaF}, we can observe an improvement in the prediction of modal coefficients incorporating the estimated eddy viscosity. Reconstructed vorticity fields are also provided in Figure~\ref{fig:2DwF} along with the variation of root mean squares error with time. Although predictions are improved with the GROM-FSM implementation compared to the GROM solution, it is observed that this improvement is only significant at later times. Moreover, some modes (especially the first few) show better results than the others, see Figure~\ref{fig:2DaF}. We believe this is caused by the assumption of fixed eddy viscosity contribution for all modes. Therefore, in the optimization process involved in FSM, higher importance is given to those first few modes (since they possess the largest contribution). Consequently, a value of eddy viscosity that yields the best correction for those first modes is computed and applied for \emph{all} modes. To mitigate this issue, we extend our closure estimation framework to allow mode-dependent variations of the eddy viscosity parameter.}

\begin{figure}[ht]
\centering
\includegraphics[width=0.95\linewidth]{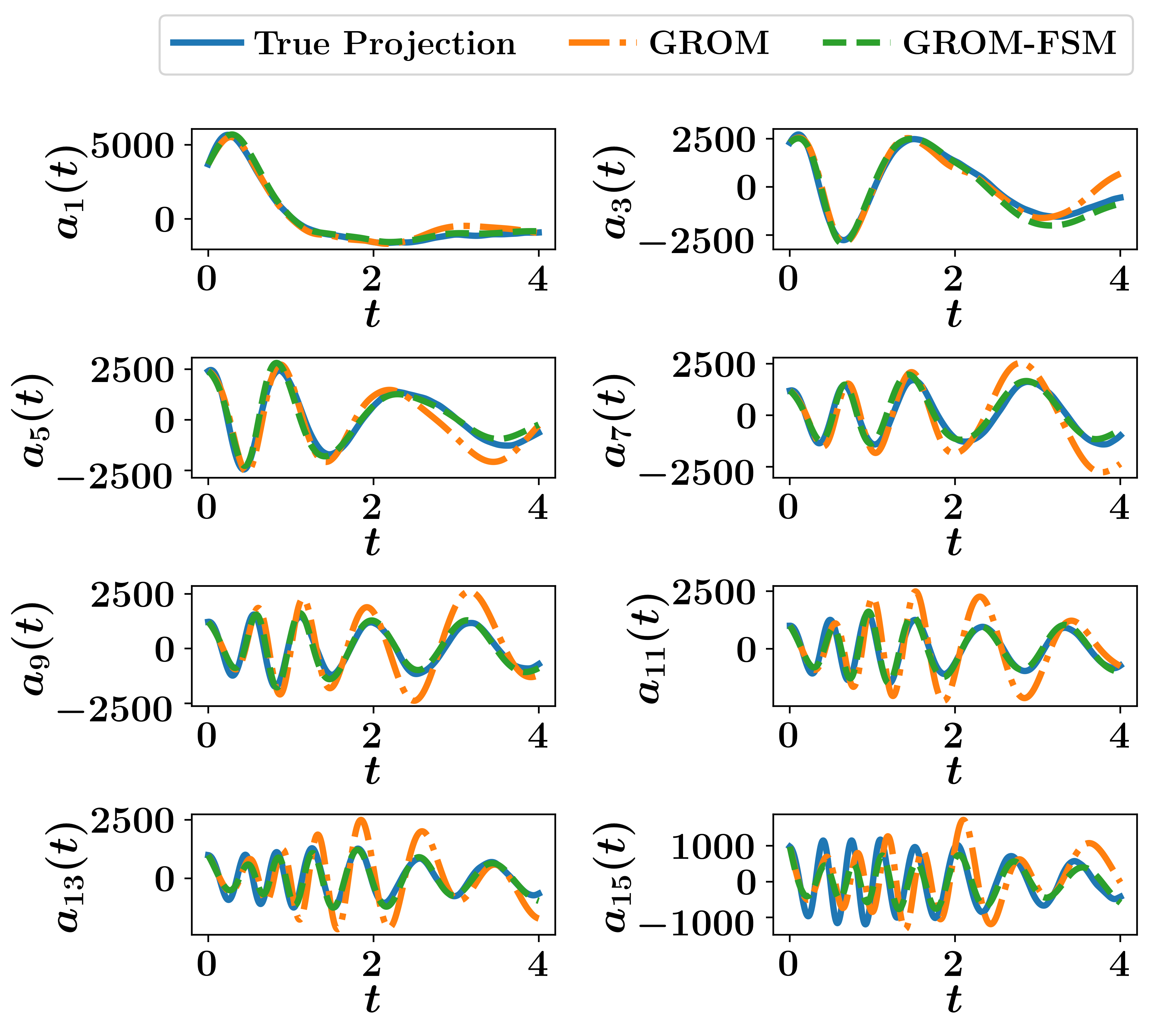}
\caption{\textcolor{rev}{Time evolution of the modal coefficients for the 2D turbulence case, with full field measurements and mode-dependent eddy viscosity closure. Note the equal improvements in predictions for all the modes.}}
\label{fig:2DaFV}
\end{figure}

\textcolor{rev}{\paragraph{\textbf{Mode-dependent eddy viscosity.}} Instead of assuming a fixed eddy viscosity parameter as is the case in Eq.~\ref{eq:gromC}, we permit the variation of this parameter with modes as follows,}
\textcolor{rev}{\begin{align}
    \dfrac{\mathrm{d}a_k}{\mathrm{d}t} &=   \mathfrak{B}_k + \nu_{e,k} \widehat{\mathfrak{B}}_k +  \sum_{i=1}^{R} \mathfrak{L}_{i,k} a_i + \nu_{e,k} \sum_{i=1}^{R} \widehat{\mathfrak{L}}_{i,k} a_i \nonumber \\
    & + \sum_{i=1}^{R} \sum_{j=1}^{R} \mathfrak{N}_{i,j,k} a_i a_j. \label{eq:gromC2}
\end{align}
Thus, our goal now is to compute the values of $\nu_{e,k}$, where $k=1,2,\dots, R$. In other words, we need to estimate $R$ local values of eddy viscosity parameters, rather than a single global value. Indeed, this approach is also common in large eddy simulations, where a spatially varying eddy viscosity is considered. Figure~\ref{fig:2DaFV} displays the time evolution of modal coefficients with a mode-dependent eddy viscosity computations. We can observe that almost equivalent improvements are obtained for all the modes, which highlights the superiority of this approach over the assumption of fixed eddy viscosity.}

\begin{figure*}[ht]
\centering
\includegraphics[width=0.95\linewidth]{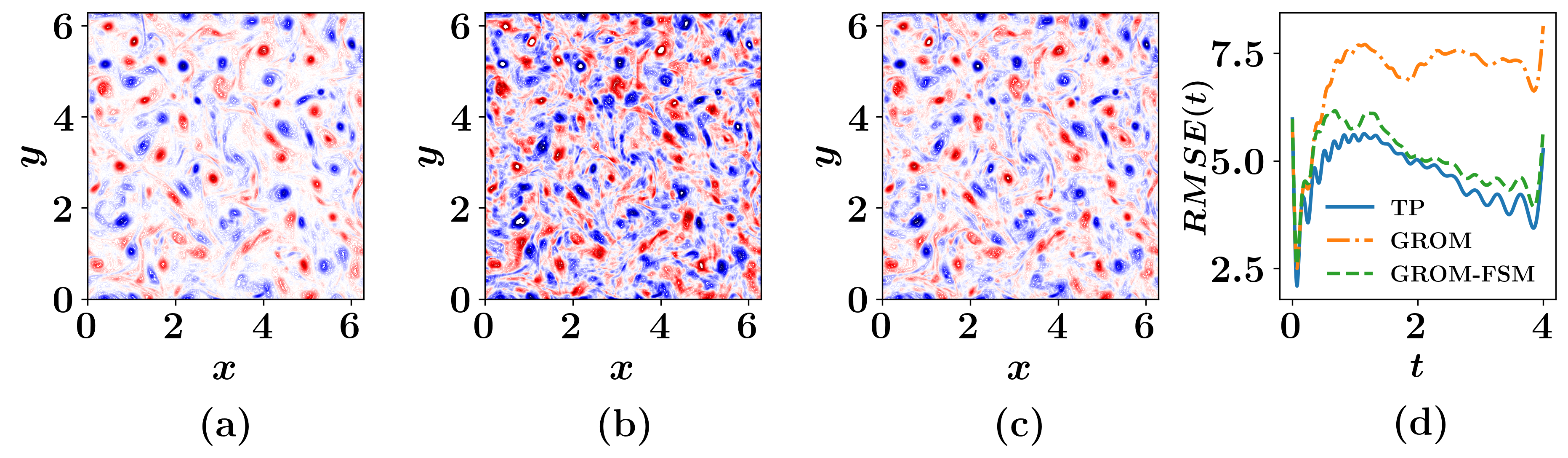}
\caption{\textcolor{rev}{Reconstructed vorticity fields for the 2D turbulence problem at $t=4$ for (a) true projection, (b) GROM, and (c) GROM-FSM along with the $RMSE$ variation with time (d), with full field measurements and mode-dependent eddy viscosity closure. Note the significant decrease in $RMSE$ at earlier times than those in Figure~\ref{fig:2DwF}.}}
\label{fig:2DwFV}
\end{figure*}

\textcolor{rev}{The quality of vorticity field reconstruction is also manifested in Figure~\ref{fig:2DwFV} with the contour plots at final time and variation of $RMSE$ with time. Interestingly, it can be noted that reductions of $RMSE$ are obtained at earlier times than those in Figure~\ref{fig:2DwF}. To understand this behavior, we investigate the GROM predictions without closure. we can see that the deviation of GROM predictions for the dynamics of the first few modes do not exhibit significant deviations during the assimilation window of $t=2$. On the other hand, the latest modes show larger deviations during the same period. However, for fixed eddy viscosity, the contribution of the first few modes (corresponding to the large convective scales) is predominant. As a result, a small value of eddy viscosity is computed to match the level of correction required for those large scales. Considering global eddy viscosity implementation, the latest modes (corresponding to small dissipating scales) receive minor corrections. On the other hand, a mode-dependent eddy viscosity implementation allows for detection of larger corrections required for disspative scales as seen in the modal coefficients predictions in Figure~\ref{fig:2DaFV} and $RMSE$ trend in Figure~\ref{fig:2DwFV}.}

\textcolor{rev}{\subsubsection{Sparse field measurement}}

\textcolor{rev}{For measurement sparsity investigation, we assume a sensor located each $32$ grid points. This corresponds to placing sensors at around $0.1\%$ of the total spatial locations. Since it has already been shown that defining a mapping from the POD coefficients to the measured field variables provides a robust data assimilation framework, we follow the same procedure here. We also consider global and local eddy viscosity implementations.}

\begin{figure}[ht]
\centering
\includegraphics[width=0.95\linewidth]{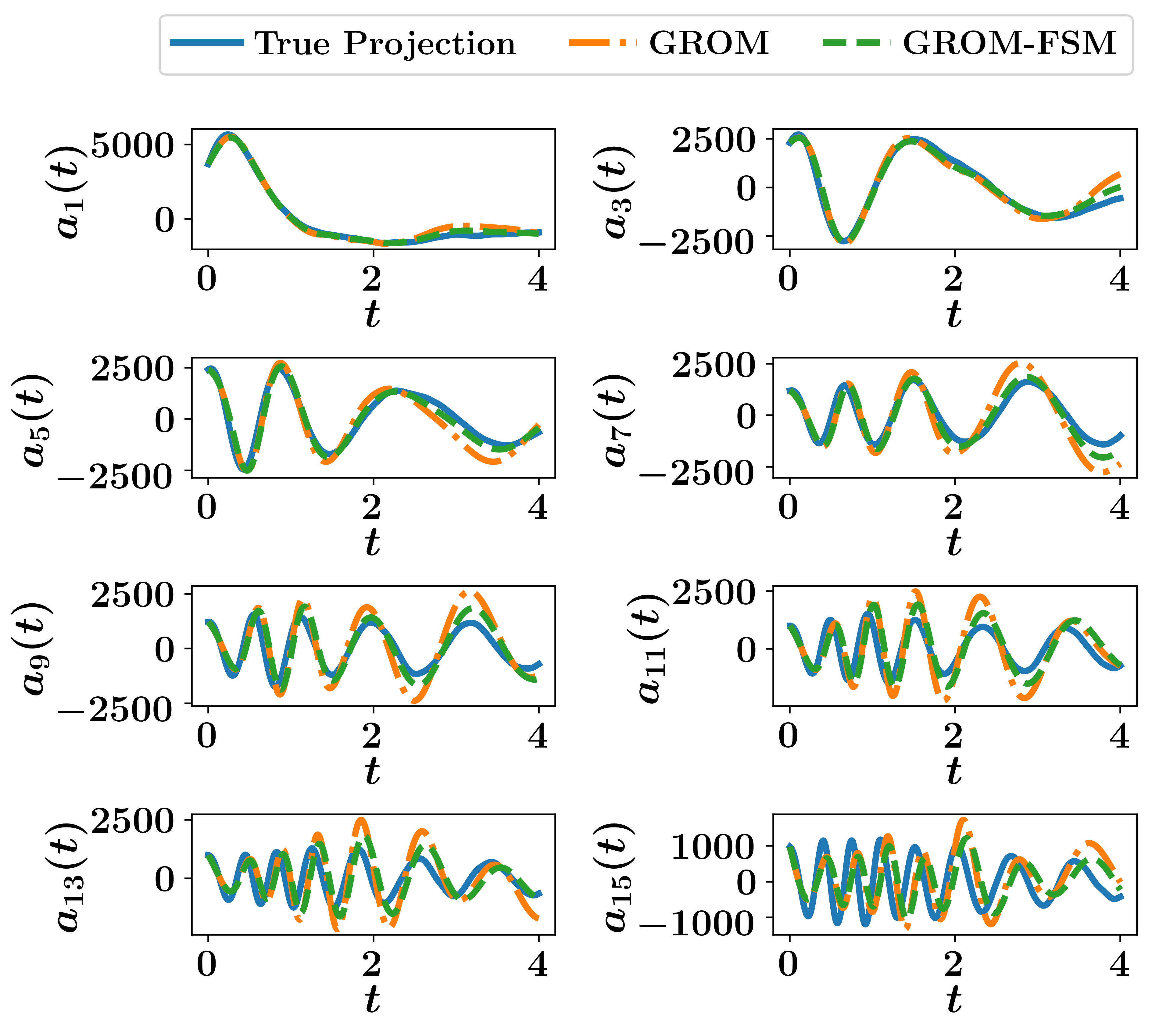}
\caption{\textcolor{rev}{Time evolution of the modal coefficients for the 2D turbulence case, with sparse field measurements and fixed eddy viscosity closure. Note the better improvements in the first few modes compared to the latest ones.}}
\label{fig:2DaS}
\end{figure}

\begin{figure*}[ht]
\centering
\includegraphics[width=0.95\linewidth]{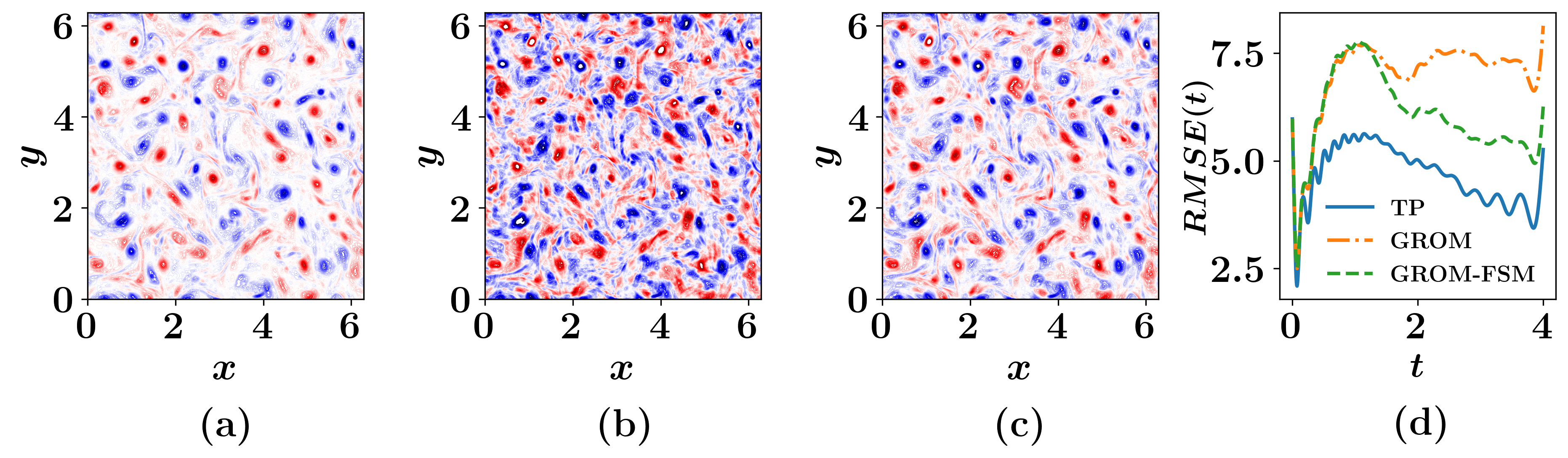}
\caption{\textcolor{rev}{Reconstructed vorticity fields for the 2D turbulence problem at $t=4$ for (a) true projection, (b) GROM, and (c) GROM-FSM along with the $RMSE$ variation with time (d), with sparse field measurements and fixed eddy viscosity closure. Reduction of $RMSE$ compared to GROM without closure starts to become remarkable around $t=2$.}}
\label{fig:2DwS}
\end{figure*}

\textcolor{rev}{\paragraph{\textbf{Fixed eddy viscosity.}} With the mapping defined in Section~\ref{sec:sp2}, we apply the FSM eddy viscosity estimation framework with sparse data and global eddy viscosity. The time dynamics of the resolved modes is demonstrated in Figure~\ref{fig:2DaS} for a few selected modal coefficients. We obtain similar results as those obtained with full field measurements, and we can observe that the predictions for the first few modes are much better than the remaining modes. Also, the reconstructed vorticity fields and computed $RMSE$ are shown in Figure~\ref{fig:2DwS}.}

\begin{figure}[ht]
\centering
\includegraphics[width=0.95\linewidth]{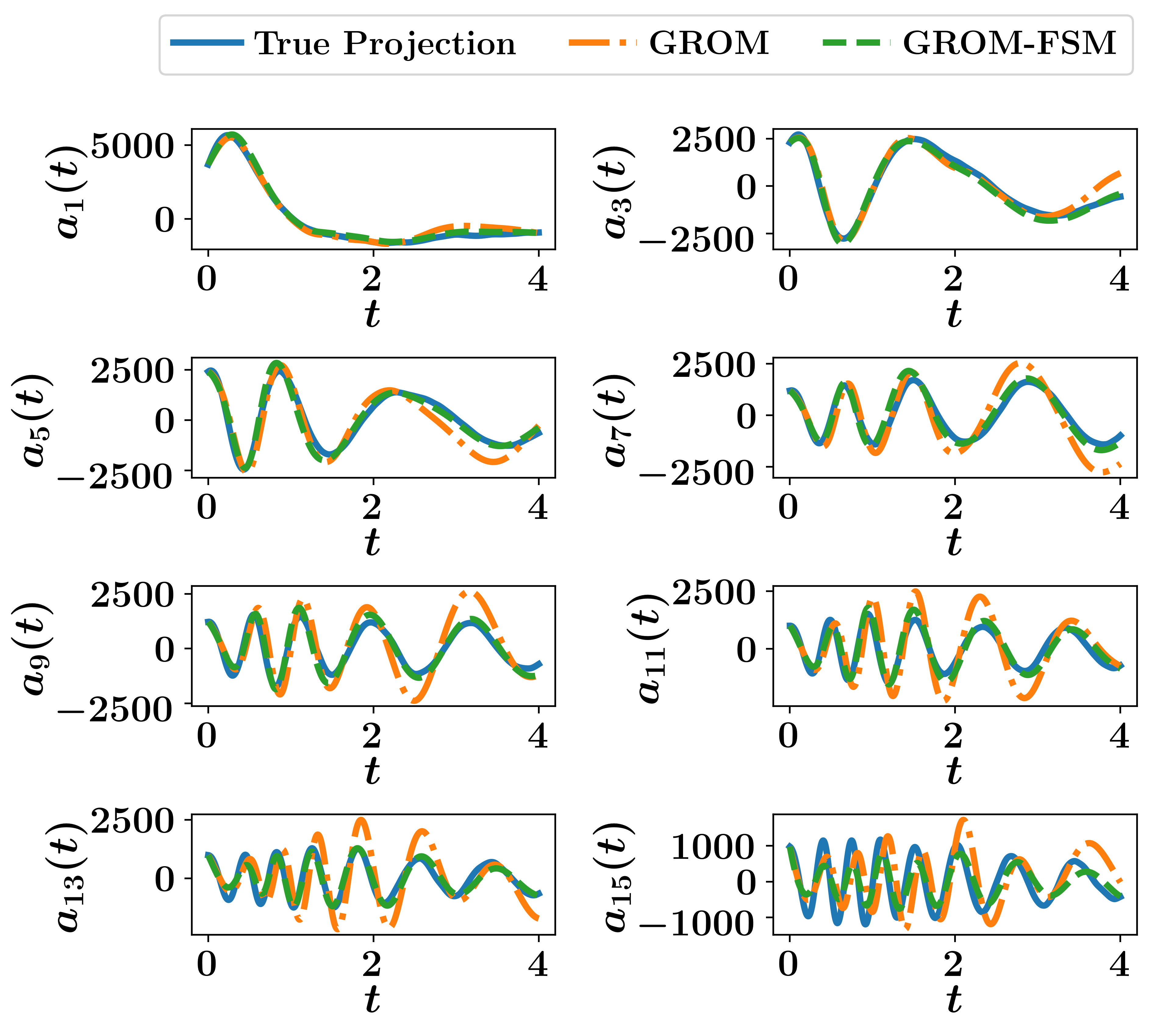}
\caption{\textcolor{rev}{Time evolution of the modal coefficients for the 2D turbulence case, with sparse field measurements and mode-dependent eddy viscosity closure.}}
\label{fig:2DaSV}
\end{figure}

\begin{figure*}[ht]
\centering
\includegraphics[width=0.95\linewidth]{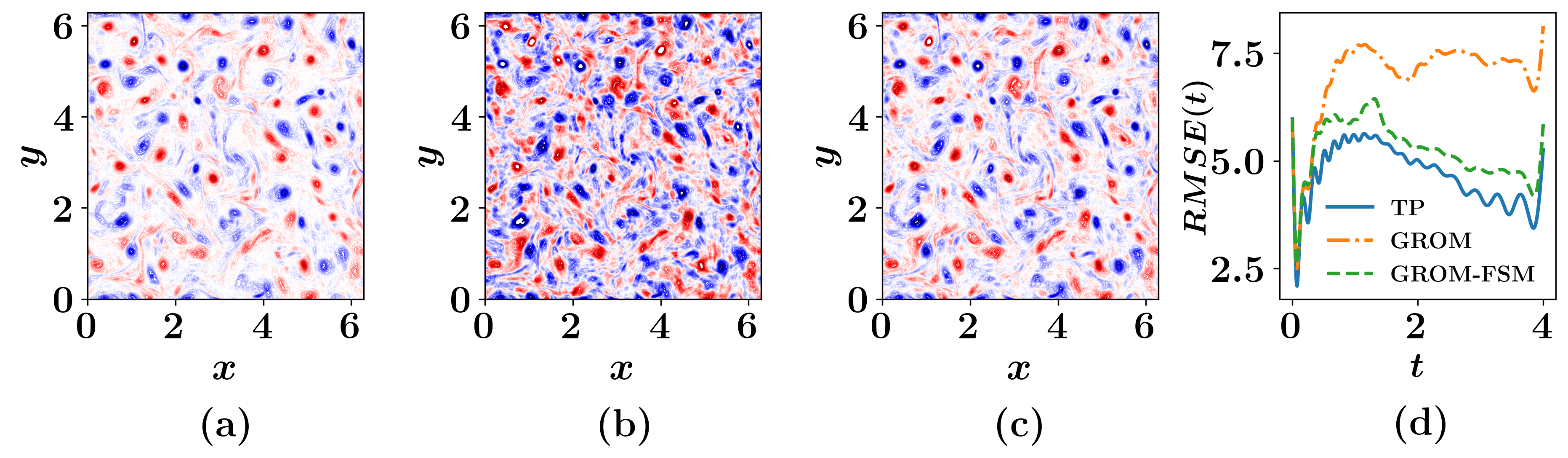}
\caption{\textcolor{rev}{Reconstructed vorticity fields for the 2D turbulence problem at $t=4$ for  (a) true projection, (b) GROM, and (c) GROM-FSM along with the $RMSE$ variation with time (d), with sparse field measurements and mode-dependent eddy viscosity closure.}}
\label{fig:2DwSV}
\end{figure*}

\textcolor{rev}{\paragraph{\textbf{Mode-dependent eddy viscosity.}} Allowing the variation of eddy viscosity yields notable enhancement of the prediction accuracy for all the modes as depicted in Figure~\ref{fig:2DaSV}, compared to Figure~\ref{fig:2DaS}. Reconstructed vorticity fields at $t=4$ with GROM, GROM-FSM and true projection results are plotted in Figure~\ref{fig:2DwSV}. We also see the reduction of $RMSE$ even with the considered $0.1\%$ sparsity in measurements data.}

\textcolor{rev}{\subsection{Computational Cost}}
\textcolor{rev}{The forward sensitivity method avoids the solution of the adjoint problem usually encountered in variational approaches for assimilating observational data to improve model's predictions. Nonetheless, the main computational burden results from the recursive matrix-matrix multiplication as described in Eq.~\ref{eq:evolutio_Vk}. However, since all computations are implemented in ROM space, the size of Jacobian matrices are $O(R)$, which reduces the memory and computational time requirements. Also, the deployed forward model in each iteration is the GROM model, which is computationally efficient when a few modes are retained in the ROM approximation. Moreover, our numerical experiments show that convergence occurs after a couple of iterations. We document the computational time for each iteration and the number of iterations for the explored test cases in Table~\ref{tab:cost1} and Table~\ref{tab:cost2} using Python implementation. We highlight here that each iteration takes around twice the time of solving the GROM equations. For instance, for the 1D Burgers problem, the solution of the GROM equations takes about $0.176 s$, while a single iteration of the FSM algorithm consumes less than $0.35 s$. Similarly, for the 2D turbulence case, the CPU time to solve the GROM equations is almost $4.445 s$, while a single iteration takes an order of $8 s$. This is comparable to the computational time of variational approaches, where the forward and adjoint problems have to be solved in each iteration. Further reductions in FSM computing time might be achieved by efficient matrix-matrix multiplication algorithms.} 

\textcolor{rev}{From Table~\ref{tab:cost2}, we can also see that the CPU time per iteration for full field measurement case is slightly smaller than that in case of sparse data. We believe that this is attributed to the difference in measurement space sizes in each case. For the full measurement case, we first project the data onto the basis functions to estimate the \emph{observed} modal coefficients, and assume our measurements live in ROM space. So, the size of observational vector and forward sensitivty matrices are all $O(R)$. However, for the sparse case where we define a map from modal coefficient to field reconstruction, the observations are assimilated in FOM space. Therefore, the size of resulting vectors and matrices, corresponding to measurements, depends on the number of sensor data. For the 2D case, this number is relatively larger than $R$, and thus the resulting matrix computations become slightly more expensive. On the other hand, for the 1D Burgers problem in Table~\ref{tab:cost1}, the CPU for either the full or sparse measurements is similar since we are using $8$ sensors for the sparse case. This is the same as the number of modes in the ROM approximation. We also notice in Table~\ref{tab:cost2} that the CPU time for mode-dependent (i.e., local) eddy viscosity estimation is larger than that in case of fixed scalar eddy viscosity approximation. This is caused by the larger sizes of the model Jacobians with respect to its parameters (see Eq.~\ref{eq:para_sen}) as well as the solution of a bigger weighted least-squares problem (e.g., Eq.~\ref{eq:FSMnue}).
}

\begin{table}[htbp!]
\caption{The CPU time (in seconds) per iteration and number of iterations required for eddy viscosity estimation using the proposed FSM-based methodology for the 1D Burgers problem. Sparse 1 refers to the implementation of mapping from measurement to POD coefficients (Section~\ref{sec:sp1}), while Sparse 2 refers to the mapping from POD coefficients to measurement (Section~\ref{sec:sp2}).}
\centering
\begin{tabular}{p{0.15\textwidth} p{0.15\textwidth} p{0.15\textwidth} }  
\hline
Measurement & CPU time (s) & No. of iterations  \\
\hline 
Full     & $0.342$ & $7$ \\ 
Sparse 1 & $0.331$ & $9$ \\ 
Sparse 2 & $0.344$ & $8$ \\ 
\hline
\end{tabular}
\label{tab:cost1}
\end{table}

\begin{table}[htbp!]
\caption{The CPU time (in seconds) per iteration and number of iterations required for eddy viscosity estimation using the proposed FSM-based methodology for the 2D turbulence problem. Here, [global] refers to the estimation of a global eddy viscosity parameter for all modes, while [local] refers to the estimation of mode-dependent eddy viscosities.}
\centering
\begin{tabular}{p{0.15\textwidth} p{0.15\textwidth} p{0.15\textwidth} }  
\hline
Measurement & CPU time (s) & No. of iterations  \\
\hline 
Full [global]    & $7.766$ & 5 \\ 
Full [local]     & $9.482$ & 4\\ 
Sparse [global]  & $7.898$ & 4 \\ 
Sparse [local]   & $9.532$ & 4\\ 
\hline
\end{tabular}
\label{tab:cost2}
\end{table}

\smallskip
\section{Concluding remarks}\label{sec:conc}
In the present study, we propose a data assimilation-based approach to provide accurate ROMs for digital twin applications. In particular, we use the forward sensitivity method (FSM) to estimate as well as update an optimal value of eddy viscosity for ROM closure. We exploit ongoing streams of observational data to improve the stability and accuracy of ROM predictions. We test the framework with the prototypical one-dimensional viscous Burgers equation characterized by strong nonlinearity \textcolor{rev}{and the two-dimensional vorticity transport equation for the 2D Kraichnan turbulence problem.} We investigate the assimilation of full field and sparse field measurements. For full field measurements, we illustrate that projecting those noisy measurements produces good estimate of \emph{observed} modal coefficients, which can therefore used to estimate an optimal value for eddy viscosity. However, we find that a similar approach of using sparse field measurements to approximate the observed states is significantly sensitive to measurements noise. On the other hand, we demonstrate that defining an observational operator via a ROM reconstruction map can be successful in utilizing sparse and noisy data. Using real-time observations can steer ROM parameters and predictions to reflect actual flow conditions. \textcolor{rev}{We also remark that the collected snapshots of full order model solutions can be assimilated by treating them as full field measurements, with negligible noise (corresponding to discretization and numerical approximation errors). This should provide a prior estimate for the eddy viscosity parameterization during an offline stage.} We emphasize that fusing ideas between physics-based closures (e.g., the ansatz for eddy viscosity) and model reduction with variational data assimilation techniques can provide valuable tools to construct reliable ROMs for long-time as well off-design predictions. This should leverage ROM implementation for real-life application.

\begin{acknowledgements}
We thank Sivaramakrishnan Lakshmivarahan for his insightful comments as well as his archival NPTEL lectures that greatly helped us in understanding the mechanics of the FSM method. This material is based upon work supported by the U.S. Department of Energy, Office of Science, Office of Advanced Scientific Computing Research under Award Number DE-SC0019290. O.S. gratefully acknowledges their support. 
Disclaimer: This report was prepared as an account of work sponsored by an agency of the United States Government. Neither the United States Government nor any agency thereof, nor any of their employees, makes any warranty, express or implied, or assumes any legal liability or responsibility for the accuracy, completeness, or usefulness of any information, apparatus, product, or process disclosed, or represents that its use would not infringe privately owned rights. Reference herein to any specific commercial product, process, or service by trade name, trademark, manufacturer, or otherwise does not necessarily constitute or imply its endorsement, recommendation, or favoring by the United States Government or any agency thereof. The views and opinions of authors expressed herein do not necessarily state or reflect those of the United States Government or any agency thereof.
\end{acknowledgements}

\bibliography{references}

\section*{Appendix A: Computing Model Jacobians}
\textcolor{rev}{We describe the computation of the model Jacobian in discrete-time formulations. We only present the case with a fixed global eddy viscosity parameter. Extension to mode-dependent closure estimation is straightforward.} For temporal discretization of the GROM equations, we use fourth-order Runge-Kutta (RK4) method as follows,
\begin{align*}
    \mathbf{a}^{k+1} &= \mathbf{a}^k + \dfrac{\Delta t}{6} (\mathbf{g}_1 + 2\mathbf{g}_2 + 2\mathbf{g}_3 + \mathbf{g}_4),
\end{align*}
where 
\begin{align*}
    \mathbf{g}_1 &= \mathbf{f}(\mathbf{a}^k, \nu_e), \\
    \mathbf{g}_2 &= \mathbf{f}(\mathbf{a}^k + \dfrac{\Delta t}{2} \cdot  \mathbf{g}_1, \nu_e), \\
    \mathbf{g}_3 &= \mathbf{f}(\mathbf{a}^k + \dfrac{\Delta t}{2} \cdot  \mathbf{g}_2, \nu_e), \\
    \mathbf{g}_4 &= \mathbf{f}(\mathbf{a}^k + \Delta t \cdot \mathbf{g}_3 ,\nu_e).
\end{align*}
Thus the discrete-time map defining the transition from time $t_k$ to time $t_{k+1}$ is written as
\begin{align*}
\mathbf{M}(\mathbf{a}^k, \nu_e) = \mathbf{a}^n + \dfrac{\Delta t}{6} (\mathbf{g}_1 + 2\mathbf{g}_2 + 2\mathbf{g}_3 + \mathbf{g}_4). 
\end{align*}
Then, the `total' Jacobian of $\mathbf{M}$ is an $R \times (R+1)$ matrix, computed as
\begin{align*}
    \mathbf{D}^k(\mathbf{M})  &= [  \mathbf{D}^k_{\mathbf{a}}(\mathbf{M}),  ~\mathbf{D}^k_{\nu_e}(\mathbf{M}) ] \\
    & = \mathbf{P} + \dfrac{\Delta t}{6} \bigg( \mathbf{D} \mathbf{g}_1 + 2 \mathbf{D}\mathbf{g}_2 + 2 \mathbf{D}\mathbf{g}_3  + \mathbf{D}\mathbf{g}_4 \bigg),
\end{align*}
where $\mathbf{P} = \begin{bmatrix} \mathbf{I}_{R}, & \mathbf{0}_{R\times1} \end{bmatrix} \in \mathbb{R}^{R\times(R+1)}$. The Jacobian of the model $\mathbf{M}$ with respect to the model state $\mathbf{a}^{k}$ is the first $R$ columns of $\mathbf{D}(\mathbf{M})$, while the Jacobian of $\mathbf{M}$ with respect to the the eddy viscosity parameter $\nu_e$ is the last column of $\mathbf{D}(\mathbf{M})$. 

Here, $\mathbf{D} \mathbf{g}_1$, $\mathbf{D}\mathbf{g}_2$, $\mathbf{D}\mathbf{g}_3$, and $\mathbf{D}\mathbf{g}_4$ are evaluated using the chain rule as follows,
\begin{align*}
    \mathbf{D}\mathbf{g}_1 &= \mathbf{D}\mathbf{f}(\mathbf{a}^k, \nu_e), \\
    \mathbf{D}\mathbf{g}_2 &= \bigg(\mathbf{D}\mathbf{f}(\mathbf{a}^k + \dfrac{\Delta t}{2} \cdot  \mathbf{g}_1, \nu_e) \bigg) \bigg(\mathbf{I}_{(R+1)}+\dfrac{\Delta t}{2}  \begin{bmatrix} \mathbf{D}\mathbf{g}_1  \\ \mathbf{Q} \end{bmatrix} \bigg), \\
    \mathbf{D}\mathbf{g}_3 &= \bigg(\mathbf{D}\mathbf{f}(\mathbf{a}^k + \dfrac{\Delta t}{2} \cdot  \mathbf{g}_2, \nu_e) \bigg) \bigg(\mathbf{I}_{(R+1)}+\dfrac{\Delta t}{2}  \begin{bmatrix} \mathbf{D}\mathbf{g}_2  \\ \mathbf{Q} \end{bmatrix} \bigg), \\
    \mathbf{D}\mathbf{g}_4 &= \bigg(\mathbf{D}\mathbf{f}(\mathbf{a}^k + \Delta t \cdot  \mathbf{g}_3, \nu_e) \bigg) \bigg(\mathbf{I}_{(R+1)}+ \Delta t \begin{bmatrix} \mathbf{D}\mathbf{g}_3  \\ \mathbf{Q} \end{bmatrix} \bigg),
\end{align*}
where $\mathbf{Q} = \mathbf{0}_{1\times(R+1)}$. Finally, the Jacobian of $\mathbf{D}\mathbf{f}(\mathbf{a}^k, \nu_e)$ is defined as  $\mathbf{D}\mathbf{f}(\mathbf{a}^k, \nu_e) = [\mathbf{D}_{\mathbf{a}}\mathbf{f}(\mathbf{a}^k, \nu_e), ~\mathbf{D}_{\nu_e} \mathbf{f}(\mathbf{a}^k, \nu_e)]$, where
\begin{align*}
    \dfrac{\partial f_k}{\partial a_j} & = (\nu + \nu_e) \mathfrak{L}_{j,k} + \sum_{i=1}^{R} \mathfrak{N}_{i,j,k} a_i + \sum_{i=1}^{R} \mathfrak{N}_{j,i,k} a_i , \\
    \dfrac{\partial f_k}{\partial \nu_e} & = \sum_{i=1}^{R} \mathfrak{L}_{i,k} a_i ,
\end{align*}
for $1 \le j,k \le R $.
 
\end{document}